\documentstyle[12pt,titlepage,newlfont]{amsart}

\newcommand{\nc}{\newcommand}

\nc{\one}{\mbox{\bf 1}}
\nc{\invtensor}{\underset{\leftarrow}{\otimes}}
\nc{\ad}{\operatorname{ad}}
\nc{\rk}{\operatorname{rk}}
\nc{\Sym}{\operatorname{Sym}}
\nc{\sym}{\operatorname{sym}}
\nc{\id}{\operatorname{id}}
\nc{\Ker}{\operatorname{Ker}}
\nc{\Aut}{\operatorname{Aut}}
\nc{\im}{\operatorname{Im}}
\nc{\ter}{\operatorname{ter}}
\nc{\intl}{\operatorname{int}}
\nc{\out}{\operatorname{out}}
\nc{\Tor}{\operatorname{Tor}}
\nc{\Hom}{\operatorname{Hom}}
\nc{\End}{\operatorname{End}}
\nc{\holim}{\operatorname{holim}}
\nc{\Ann}{\operatorname{\bf Ann}}
\nc{\lwt}{\operatorname{lwt}}
\nc{\rwt}{\operatorname{rwt}}
\nc{\wt}{\operatorname{wt}}
\nc{\ch}{\operatorname{ch}}
\nc{\Ens}{{\cal E}\mbox{\em ns}}
\nc{\Sch}{{\cal S}\mbox{\em ch}}
\nc{\Tot}{\operatorname{Tot}}
\nc{\Th}{\operatorname{Th}}
\nc{\Cech}{\check{C}}
\nc{\Spec}{\operatorname{Spec}}
\nc{\Prim}{\operatorname{Prim}}
\nc{\Fract}{\operatorname{Fract}}

\nc{\xa}{y^{-1}_{-\alpha}}
\nc{\xb}{y^{-1}_{-\alpha}y_{-\alpha-\beta}}
\nc{\xc}{\ol{y}_{-\alpha-\beta}}
\nc{\xd}{y_{-\gamma}}
\nc{\xe}{y^{-1}_{-\alpha}y_{-\alpha-\beta-\gamma}}
\nc{\xf}{\ol{y}_{-\alpha-\beta-\gamma}}
\nc{\xg}{y_{-\beta}}
\nc{\xh}{y_{-\beta-\gamma}}
\nc{\xy}{y_{-\alpha-2\beta-\gamma}}
\nc{\xj}{\ol{y}_{-\beta-\gamma}}

\nc{\Dglie}{\operatorname{{\cal D}glie}}
\nc{\Dgcoalg}{\operatorname{{\cal D}gcoalg}}
\nc{\Homcoalg}{\operatorname{{\cal H}omcoalg}}
\nc{\Homlie}{\operatorname{{\cal H}omlie}}
\nc{\Dgmod}{\operatorname{{\cal D}gmod}}
\nc{\Dgcomod}{\operatorname{{\cal D}gcomod}}

\nc{\pa}{\partial}
\nc{\co}{\cal O}

\nc{\fg}{\frak g}
\nc{\fn}{\frak n}
\nc{\fb}{\frak b}
\nc{\CO}{\cal O}
\nc{\fe}{{\frak n}_+}
\nc{\fh}{\frak h}
\nc{\ft}{\frak t}

\nc{\dirlim}{\underset{\rightarrow}{\lim}\,} 
\nc{\nen}{\newenvironment}
\nc{\ol}{\overline}
\nc{\ul}{\underline}
\nc{\ra}{\rightarrow}
\nc{\lra}{\longrightarrow}
\nc{\Lra}{\Longrightarrow}
\nc{\Lla}{\Longleftarrow}
\nc{\Llra}{\Longleftrightarrow}
\nc{\hra}{\hookrightarrow}
\nc{\iso}{\overset{\sim}{\lra}}

\nc{\Thm}[1]{Theorem~\ref{#1}}
\nc{\Prop}[1]{Proposition~\ref{#1}}
\nc{\Lem}[1]{Lemma~\ref{#1}}
\nc{\Cor}[1]{Corollary~\ref{#1}}
\nc{\Conj}[1]{Conjecture~\ref{#1}}
\nc{\Claim}[1]{Claim~\ref{#1}}
\nc{\Defn}[1]{Definition~\ref{#1}}
\nc{\Exa}[1]{Example~\ref{#1}}
\nc{\Rem}[1]{Remark~\ref{#1}}
\nc{\Note}[1]{Note~\ref{#1}}

\nen{thm}[1]{\label{#1}{\bf Theorem.\ } \em}{}
\nen{prop}[1]{\label{#1}{\bf Proposition.\ } \em}{}
\nen{lem}[1]{\label{#1}{\bf Lemma.\ } \em}{}
\nen{cor}[1]{\label{#1}{\bf Corollary.\ } \em}{}
\nen{conj}[1]{\label{#1}{\bf Conjecture.\ } \em}{}
\nen{claim}[1]{\label{#1}{\bf Claim.\ } \em}{}

\nen{defn}[1]{\label{#1}{\bf Definition.\ } }{}
\nen{exa}[1]{\label{#1}{\bf Example.\ } }{}

\nen{rem}[1]{\label{#1}{\em Remark.\ } }{}
\nen{note}[1]{\label{#1}{\em Note.\ } }{}
\nen{exer}[1]{\label{#1}{\em Exercise.\ } }{}

\begin{document}
\title[]{The prime spectrum
of a quantum Bruhat cell translate}
\author[]{Maria Gorelik}
\address{ Dept. of Theoretical Mathematics,
 The Weizmann Institute of Science, 
Rehovot 76100, Israel,
{\tt email: remy@@wisdom.weizmann.ac.il} 
}
\begin{abstract}
The prime spectra of two families of algebras,
$S^w$ and $\check{S}^w$, $w\in W,$ 
indexed by the Weyl group $W$ of a semisimple finitely dimensional
Lie algebra $\frak g$, are studied in the spirit of~\cite{j1}.
The algebras $S^w$ have been introduced by A.~Joseph (see~\cite{j4}, 
Sect.~3). They are $q$-analogues of the algebras of regular functions 
on $w$-translates of the open Bruhat cell of a semisimple Lie group $G$ 
corresponding to the Lie algebra $\fg$. 

We define a stratification of the spectra into
components indexed by pairs $(y_1,y_2)$ of elements of the Weyl group
satisfying $y_1\leq w\leq y_2$. Each component admits a unique minimal
 ideal
which is explicitly described. We show the inclusion relation
of closures to be that induced by Bruhat order.
\end{abstract}

\thanks{The work was partially supported by the Hirsch and Braine 
Raskin Foundation}

\maketitle

\newpage
\begin{center}
Running head: QUANTUM BRUHAT CELL TRANSLATE
\end{center}
\newpage

\section{Introduction}

In this work we study the prime spectra of two families of algebras,
$S^w$ and $\check{S}^w$, $w\in W,$ 
indexed by the Weyl group $W$ of a semisimple finitely dimensional
Lie algebra $\frak g$.
The algebras $S^w$ have been introduced by A.~Joseph (see~\cite{j4}, 
Sect.~3). They are $q$-analogues of the algebras of regular functions 
on $w$-translates of the open Bruhat cell of a semisimple Lie group $G$ 
corresponding to the Lie algebra $\fg$. 

The corresponding classical objects, the algebras of regular functions on 
different $w$-translates of the open Bruhat cell, are isomorphic 
to each other polynomial algebras of rank $\dim\fn^+$.
  
The $q$-analogues $S^w$ are much more interesting. For instance,
their centres have different Gelfand-Kirillov dimension for different $w\in W$ 
--- see~\Rem{nonisomorphic}. In particular, $S^w$ are not in general
isomorphic for different $w\in W$.

The algebras $S^w$ admit a structure of right $U_q(\fg)$ module which comes
from the right action of $U_q(\fg)$ on the quantum function ring $R_q[G]$.
The action of the root torus $T\subseteq U_q(\fg)$ on $S^w$ can be naturally
extended to an action of the weight torus $\check{T}\supseteq T$. The second
family of algebras, $\check{S}^w$, are obtained as the skew-products
$ \check{S}^w=S^w\#\check{T}.$

The starting point of the construction of the rings $S^w$ is the ring $R^+$
which is a quantization of the ring of global regular functions on the 
``base affine space'' $G/N$, see~\cite{j4}, 1.2. The algebra $S^w$ is
obtained as a zero weight space of a localization of $R^+$. This is why 
the rings $S^w$, $ \check{S}^w$  are denoted almost everywhere as $R^w_0$,
$\check{R}^w_0$ respectively.

In the case $w=e$ the algebra $S^e$ is isomorphic to the quantized
enveloping algebra $U_q(\fn^-)$ of the maximal nilpotent subalgebra 
$\fn^-\subseteq\fg$ --- see~\cite{j4}, 3.4. The corresponding skew-product
algebra $\check{S}^e$ is isomorphic to $\check{U}_q(\fb^-)$.

The prime spectrum of the algebra $\check{S}^e\cong\check{U}_q(\fb^-)$
was described by A.~Joseph ~\cite{j1}, Sect.9. It is presented as
a disjoint union of locally closed strata  $X(w)$ indexed by the
elements of the Weyl group. Moreover, the strata $X(w)$ admit an action
of a group ${\Bbb Z}_2^l\subseteq\Aut(\check{S}^e)$ and the quotient
$X(w)/{\Bbb Z}_2^l$ is isomorphic (as a partially ordered set) to the
spectrum of a Laurent polynomial ring.

In this paper we present a similar description
(\Prop{dcms}) of the spectrum of $\check{S}^w$ for arbitrary $w\in W.$ 
In our case the strata $X_w(y,z)$ are indexed by a more complex set: this is 
the collection 
$$W\overset{w}{\diamond}W:=\{(y,z)\in W\times W|\ y\leq w\leq z\}$$
where $\ \leq\, $ is the Bruhat order. Note that $\, W\overset{w}{\diamond}W$ 
inherits an order relation through 
$\ (y,z)\succeq (y',z') \text{ iff } y\leq y',\ z\geq z'.$
In~\Cor{stratorder} we prove that the closure of $X_w(y,z)$ coincides with 
the union
of $X_w(y',z'):\ (y,z)\succeq (y',z')$.

The spectrum of $S^w$ is a union of strata $Y_w(y,z)$ indexed by the same
set $W\overset{w}{\diamond}W$ (\Prop{dcms}). One has also a similar
decomposition of a ``generic part'' $\Spec_+R^+$ of the spectrum of $R^+$
(see~\ref{specplus},\Cor{dcmrplus}). Here the strata $X(y,z)$ are indexed
by the set $$W{\diamond}W:=\{(y,z)\in W\times W|\ y\leq z\}.$$
The strata $X_w(y,z)$ (resp., $Y_w(y,z)$) are isomorphic for different 
$w:\ y\leq w\leq z$ (\Prop{prpisoxx}).
Moreover, $X_w(y,z)$ are all isomorphic to the component $X(y,z)$
of $\Spec_+R^+$ (\Prop{prpisoxy}). It turns out that the component $X(y,z)$ 
is isomorphic
(up to an action of a group ${\Bbb Z}_2^l$) to the spectrum of a Laurent
polynomial ring --- see~\Thm{spectra}.

The stratum $X(y,z)$ admits a unique minimal 
element $Q(y,z)$ which we calculate explicitly in~\Prop{mini}.
We deduce from this that
the stratum $Y_w(y,z)$ also admits a unique minimal element $Q(y,z)_w$
which can be expressed through a localization of $Q(y,z)$ (\Cor{minx}). 
Then the unique
minimal element of the stratum $X_w(y,z)$ can be written as
$Q(y,z)_w\#\check{T}$ --- see~\Cor{minx}. 
The prime ideals $Q(y,z),\ Q(y,z)_w,\ Q(y,z)_w\#\check{T}$ are completely
prime.

In the last Section~\ref{ccc} we calculate the centres of the rings $S^w$
(note that the centres of $\check{S}^w$ are trivial). These are polynomial
rings whose dimension depends on $w\in W$. 

In the special case $\fg=\frak{sl}_4$ the prime and the primitive
spectra of $S^w$ were calculated in~\cite{g1}. The results of the
first draft of this paper have been announced in~\cite{g2}.

{\em Acknowledgement.} I am greatly indebted to Prof. A.~Joseph who posed 
the problem. His book "Quantum groups and their primitive ideals" was the 
main inspiration of the present work. I am also grateful to him for 
reading of the first draft of the manuscript and for numerous suggestions.
I am grateful to V.~Hinich for helpful discussions and support.

\section{The rings $S^w,\ \check {S}^w$}
\subsection{}
\label{ch} 
The base field $k$ is assumed to be of characteristic zero and $K$ is an
extension of $k(q)$.
 
Let $\frak g$ be a semisimple Lie algebra and $\, U_q(\frak g)\,$ be the 
Drinfeld-Jimbo
quantization of $\, U(\frak g)\,$ defined for example in~\cite{j}, 3.2.9 whose
notation we retain. In this 
$U_q(\frak g)\;$ is a $K$-algebra generated by $x_i$, $y_i$, $t_i$, $t_i^{-1}$ 
$i=1,\ldots ,l\;$ where $l$ is the rank of $\frak g$. 
Denote the extension of $\, U_q(\frak g)\,$ by the torus 
$\check{T}$ of weights (\cite{j}, 3.2.10) 
by $\check {U_q}(\frak g)\;$. Consider the subalgebra $U_q(\frak n^-)\;$
generated by the $y_i$, $i=1,\ldots ,l\;$ (\cite{j}, 3.2.10). 
By~\cite{j}, 10.4.9 $U_q(\frak n^-)\;$ admits a structure of a right 
$U_q(\frak g)$-module such that:

(1) This module structure is compatible with the algebra 
structure of $U_q(\frak n^-)\;$ and the coproduct on $\, U_q(\frak g)\,$.
  
(2) Endowed with this $U_q(\frak g)$-module  
structure $U_q(\frak n^-)\;$ is isomorphic to the dual $\delta M(0)\;$ 
of the  $U_q(\frak g)$-module Verma (~\cite{j}, 5.3) of highest weight zero.  

After Lusztig-Soibelman the braid group of $\frak g\;$ acts on  
$U_q(\frak g)\;$ 
by automorphisms $r_w\;$ such that if $\tau (\lambda )\;$ is an element of the
torus $T\;$ and $\overline w\ $ is the image of $w\ $ in the Weyl group $W\;$ 
of $\frak g\;$ then:
$$r_{w} \tau (\lambda )=\tau (\overline w \lambda ).$$ 
Fix an element $\overline w\ $ of the Weyl group and let $w\ $ be a 
representative
of $\overline w\ $ in the braid group. The automorphism $r_w\;$ acts on 
the category of $U_q(\frak g)$-modules by transport of structure. Denote 
$\;\left (\delta M(0)\right )^{r_w}$ by $S^w$. As noted
in~\cite{j}, 10.4.9 the $\check{T}$-character of $S^w$ is given by the 
formula  
$$\ch S^w=w\left(\prod_{\beta\in\Delta^-}(1-e^{\beta})^{-1}\right)
=\prod_{\beta\in w\Delta^-}(1-e^{\beta})^{-1}.$$
Suppose $\psi$ is an automorphism of $U_q(\frak g)$ such that the module
$\left (\delta M(0)\right )^{\psi}$ has the same character as $S^w$.  
Then the module $N=\left (\delta M(0)\right )^{r_w^{-1} \psi }$ has the same
character as $\delta M(0)$. Since $N$ is obtained from $\delta M(0)$ by 
transport
of structure the following property of $\delta M(0)$ holds also for $N$: 
if $v_0$
is a vector of weight zero and $v$ is a vector of $N$ then $v_0$ belongs to the
submodule generated by $v$. Hence the dual module $\delta N$ is 
generated by a highest weight vector. Yet it is also has the same character as 
the Verma module $M(0)$, so $\delta N$ is isomorphic to $M(0)$, $N$ is 
isomorphic
to $\delta M(0)$ and $\left (\delta M(0)\right )^{\psi }\;$ is isomorphic
to $\left (\delta M(0)\right )^{r_w}$. Hence the $U_q(\frak g)$-module $S^w$ 
depends 
only on the class $\overline w\ $ of $w\ $ in the Weyl group $W$ of $\frak g$.
   
According to ~\cite{j}, 10.2.9, $S^w$ admits the structure of a 
$U_q(\frak g)$-algebra and this further extends to a 
$\check {U_q}(\frak g)\ $-algebra structure. Moreover one checks that
the $\check {U_q}(\frak g)\ $-algebra structure 
on the module $S^w$ is uniquely determined up to a scalar by its module 
structure
and the requirement that a non-zero vector of weight zero is the identity
of the ring (see also~\cite{k}, prop. 3.2). 
The automorphism $r_w$ is an algebra automorphism but it does not preserve the
coalgebra structure of $U_q(\frak g)$. Thus one should not expect that
the algebras $S^w\;$ are isomorphic for different elements 
$\overline w\in W\,$. Rather we obtain a collection of $U_q(\frak g)$-algebras 
parametrized by $W$ which are generally distinct. Trying to 
understand the possible isomorphisms between them was a main
motivation for our present work. Our results suggest that $S^w$ is
isomorphic to $S^{w'}\,$ iff $W\overset{w}{\diamond}W$ and
$W\overset{w'}{\diamond}W\,$ are isomorphic as
ordered sets.

\subsection{}
\label{dcpm}
Let $w_0$ be the longest element of the Weyl group. 
Consider the involution $\psi$ of the algebra $U_q(\frak g)$ defined by
$$\psi (x_i)=-y_i\ \ \ \ \ \ \psi (t_i)=t_i^{-1}.$$
Then by the character formula of~\ref{ch} one has 
$$\ch \left ( S^w\right )^{\psi }=\ch S^{w w_0}.$$  
By the reasoning of~\ref{ch} the modules $\left ( S^w\right )^{\psi }$ and 
$S^{w w_0}$ are isomorphic and hence are isomorphic as algebras. 
The map $\psi\;$ is an algebra automorphism and coalgebra 
antiautomorphism. The last implies that the  $U_q(\frak g)$-algebras $S^w$ and 
$\left ( S^w\right )^{\psi }$ have opposite algebra structures.
Hence the algebras $S^w$ and $S^{w w_0}$ are opposites.

\subsection{}
\label{second}
Fix a triangular decomposition 
$\, \frak g=\frak n^-\oplus\frak h\oplus\frak n^+\,$
and let $\,\pi=\{\alpha_1, \alpha_2, \dots ,\alpha_l\}\ $ be 
the corresponding set of simple roots. Let $\, Q(\pi)=\Bbb Z\pi\,$,
$Q^{\pm}(\pi)=\pm\Bbb N\pi\,$, $P(\pi)$ (resp., $P^+(\pi)$) be the set of 
weights 
(resp., dominant weights) and $\,\{\omega _i\}_ {i=1}^{l}\ $ be the set of
fundamental weights. Define an order relation on $P(\pi)\,$ 
by $\,\mu \geq\nu\,$ if $\, \mu-\nu\in Q^+(\pi).$
Let $\tau$ be the isomorphism of the additive
group $Q(\pi)$ to themultiplicative group $T$ defined by 
$\tau(\alpha_i)=t_i,\ i=1,\ldots,l$.
We can extend $\tau$ to the isomorphism of $P(\pi)$ onto $\check {T}$.

For each $\,\lambda\in P^+(\pi)\,$ let $\, V(\lambda)\,$ be
the $U_q (\frak g)$ module with highest weight $\lambda$ and
$\,c^\lambda_{\xi, v}\,$: $\xi\in  V(\lambda)^*$, $v\in V(\lambda)\,$ be the 
element
$a\mapsto\xi (av)$ of $U_q(\frak g)^*\,$. Let $R_q [G]\,$ be the Hopf 
subalgebra of 
$U_q(\frak g)^*\,$ generated as a vector space by these elements. 
By~\cite{j}, 9.1.1 $R_q[G]$ 
admits a structure of a $U_q (\frak g)$-bialgebra.

Let $u_\lambda$ be a highest weight vector of 
$V(\lambda)$ and $V^+(\lambda)\,$
denote the subspace of $R_q[G]$ generated by the 
$\,c^\lambda_{\xi, u_\lambda}\,$:
$\xi\in  V(\lambda)^*\,$. Then 
$\, R^+:=\oplus_{\lambda\in P^+ (\pi)} V^+(\lambda)\,$
is a subalgebra of $R_q[G]\,$. Moreover $R^+$ is a right
$U_q (\frak g)$-submodule and left $T$-submodule of  $R_q [G]\,$. The left
$T$-action defines a $P^+(\pi)$-grading on $R^+$. Indeed the weight subspace 
of weight $\, \lambda$
is just $V^+(\lambda)\,.$ Hence $V^+(\lambda)\,$ is invariant 
with respect to the right action of $U_q (\frak g )$ and the 
multiplication satisfies the Cartan multiplication rule:
$$V^+(\mu) V^+(\lambda)= V^+(\lambda +\mu).$$

Let $\ \Omega ( V^+(\lambda))\ $ denote the set of weights of $V^+(\lambda)\, $
for the right $T$-action counted with their multiplicites. (This is just the 
set of weights of  $V(\lambda)$).

For each $\,w\in W$ let $\xi_{w\lambda}\,$ be a 
vector of the weight $w\lambda\,$ in $V(\lambda)^*\,$ viewed as a right 
$U_q (\frak g)$ module and write $\,c^\lambda_{\xi_{w\lambda}, u_\lambda}\,$
 (resp., 
$\,c^\lambda_{\xi, u_\lambda}\,$) simply as $c^\lambda _w\,$ (resp., 
$c^\lambda _\xi\,$). The elements $c^\lambda _w\,$ are defined up to scalars.
By~\cite{j}, 9.1.10 these scalars can be chosen so 
that $\,c^{\mu} _w c^{\nu} _w=c^{\mu +\nu}_w$
for any $\,\mu ,\nu\in P^+ (\pi)\ $ and 
$c_w=\{c^\lambda _w : \lambda\in P^+(\pi)\}$
becomes an Ore set in $R^+\,$. Extend $\,c^{\mu} _w$ to $\mu\in P(\pi)\,$
through 
$\ c^{\mu -\nu}_w=c^{\mu} _w (c^\nu _w)^{-1}\ \ \forall \mu ,\nu\in P^+ (\pi)$.

Consider the localized algebra $\, R^w:=R^+[c_w^{-1}]$; 
by~\cite{j}, 4.3.12 the right action of $U_q (\frak g )$ extends to $R^w$.
Since each of $c^\lambda _w$ is homogeneous it follows that the 
$P^+(\pi)$-grading
on $R^+$ extends to a $P(\pi)$-grading on $R^w$; again the homogeneous 
components
are invariant  with respect to the right action of $U_q (\frak g)$.
It implies that the zero weight subspace $R^w_0\,$ of $R^w$
with repsect to the left action of $T$ is a $U_q (\frak g)$-subalgebra
of $R^w\,$ and as suggested in~\cite{j4}, 3.1, it may be viewed as 
a $q$-analogue of the
algebra of regular functions on the $w$-translate of the open Bruhat cell.
Since $R^+\,$ is a domain of finite Gelfand-Kirillov dimension it admits
a skew-field of fractions and this contains the $\, R^w: w\in W.$ 
Again $\ c^{-\lambda} _w V^+(\lambda)\hookrightarrow
 c^{-(\lambda+\nu)} _w V^+(\lambda+\nu)\ \  \forall \lambda, \nu \in P^+(\pi).$
Thus one may write
\begin{equation}
R_0^w=\sum_{\lambda\in P^+(\pi)} c^{-\lambda}_w V^+(\lambda)\cong
\underset{\lambda\in P^+(\pi)}{\dirlim} c^{-\lambda}_w V^+(\lambda).
\label{A}
\end{equation}
This implies that the rings of fractions of $R_0^w$ coincide for different $w.$

By~\cite{j}, 10.4.8 $S^w\,$ and $R^w_0\,$ are isomorphic as a 
$U_q(\frak g)$-algebras.

Denote by $\check {R}^w_0$ the skew-product of $R^w_0$ and the fundamental 
torus
$\check {T}$ through the action of $\check {T}$ on 
$\check {U_q}(\frak g)$-module $R^w_0$--- see~\ref{defskewpr}.

\subsection{}
\label{noeth}
By~\cite{j4}, 6.4, 6.6, $R^+\,$ and $S^w$ are left and right noetherian. 
By~\cite{mr}, 2.9 it follows that $\check {R}^w_0$ is also noetherian.

\subsection{}
\label{a3}
Set $w=e$. Then $S^e\;$ is isomorphic to $U_q(\frak n^-)$ as a 
$\check {U_q}(\frak g)$-algebra.
Consider the subalgebra $\check {U_q}(\frak b ^-)\;$ of 
$\check {U_q}(\frak g)\;$ which is the skew-product of $U_q(\frak n^-)\;$  and 
the fundamental torus $\check{T}.$ The algebra $\check {U_q}(\frak b ^-)\;$
can be also considered as the skew-product of $S^e$ and $\check{T}\,$
through the action of $\check{T}\,$ on 
$\check {U_q}(\frak g)$-module $S^e$. By~\cite{j}, 10.1.11 it follows that
the isomorphism~\ref{second} of $S^e\simeq U_q(\frak n^-)$ with 
$R^e_0$ extends to an isomorphism of $\check {U_q}(\frak b ^-)$  with $R^e.$

By~\cite{j1}, Sect.10 the prime and
primitive spectra of $\check {U_q}(\frak b ^-)\;$ take the following form
$$\Spec  \check {U_q}(\frak b ^-)=\coprod_{w\in W} X(w)\ ,$$  
$$\Prim   \check {U_q}(\frak b ^-)=\coprod_{w\in W} X^{max}(w)\ ,$$  
where each $X(w)$ is the spectrum of some Laurent polynomial ring up to an
action of $\Bbb{Z}_2^l$ and all prime ideals are completely prime. 

Each $X(w)$ has a unique minimal element $Q(w)$
which has the following nice description in the notation of~\ref{second}.
Fix $w\in W$. For each $\,\lambda\in P^+(\pi)\,$ let 
$u_{w\lambda}\in V(\lambda)\,$ be
a vector of the weight $w\lambda\,$. Denote by $V_w^+(\lambda)^{\perp}\,$ the
orthogonal of the Demazure module 
$V_w^+(\lambda):=U_q(\frak b ^+)u_{w\lambda}\,$ 
in $V(\lambda)^*\,$, the latter identified with $V^+(\lambda)\,$. 
Then~\cite{j}, 10.1.8 
$$Q(w)=\sum_{\lambda\in P^+(\pi)} V_w^+(\lambda)^{\perp}.$$

\subsection{}
\label{normal}
An element $x$ of a ring $A$ is called {\em normal} if $xA=Ax$. If $A$ is 
prime a non-zero normal element is regular. Each regular normal
element determines an automorphism of the ring sending $a\in A\,$
to the unique element $b\in A\,$ such that $\, xa=bx$. 

Let $A$ be a ring, $x$ be an element of $A$ and $c$ be a subset of $A$.
Suppose that the multiplicative closures of $c$ and $\left\{ x \right\}$ 
are Ore sets in $A$. In this case we denote the localizations of the ring $A$ 
at the corresponding multiplicative closures respectively by $A[c^{-1}]$, 
$A[x^{-1}]$.

\section{two lemmas}
\subsection{}
Let $S$ be an algebra graded by a free abelian group $H$. Then

\begin{lem}{l2}
(i) A graded ideal $P$ is prime iff
for any homogeneous $a,b\in S\backslash P$ there exists
$c$ such that $acb\not\in P$. 

(ii) Take a prime ideal $I$ of $S$ and let $J$ be
a maximal homogeneous ideal contained in $I$.
Then $J$ is prime. 
\end{lem}
\begin{pf}
(i) Assume that for any homogeneous $a,b\in S\backslash P$ there exists
$c$ such that $acb\not\in P$. Take any $a',b'\in S\backslash P$.
We can assume that none of the homogeneous 
components of $a'$ and of $b'$ belong to $I$. 
Fix a lexicographic order on $H$.
Denote by $a$ (resp., $b$) the minimal homogeneous 
component of $a'$ (resp., $b'$) with respect to the order. 
Take $c$ such that $acb\not\in P$.
Then $ac'b\not\in P$ for some homogeneous component $c'$ of $c$.
Since the minimal homogeneous component of $a'c'b'$ is just $a c' b$, 
it follows that $a'c'b'\not\in P$ so $P$ is prime as required. 

(ii) Observe that $J$ is a linear span of the set of homogeneous
elements of $I$. Take homogeneous $a,b\not\in J$. Then
$a,b\not\in I$ so there exists $c$ such that $acb\not\in I$. 
Hence $acb\not\in J$ that, by (i), gives the required assertion. 
\end{pf} 

\subsection{}
\label{defskewpr}
Let $S$ be a $K$-algebra, $\check {T}\cong \Bbb Z ^l$ be a torus acting on $S$
by right automorphisms. Denote the action of $t\in \check {T}$ on $s\in S$
by $s.t$. Define an algebra structure on $S\otimes K[\check {T}]$ through
$$(s_1\otimes t_1)(s_2\otimes t_2)=(s_1(s_2.t_1^{-1})\otimes t_1t_2).$$
The vector space $S\otimes K[\check T]$ endowed with the above
algebra structure is called the skew-product $S\#\check {T}$. It
will be denoted also by $\check {S}$.
Denote by $\left(\Spec S \right)^{\check {T}}$ the set of
$\check{T}$-invariant prime ideals of $S$.
 
\begin{lem}{l3}
(i) If $I\in \left(\Spec S\right)^{\check {T}}\,$ then
$J:=(I \#\check{T})$ is prime in $\check {S}\ $ and $\, J\cap S=I.$ 

(ii) Assume that  $\check {T}$ acts on $S$ by semisimple automorphisms
and the set of weights $H$ is a subset of a free abelian group.
If $J\in\Spec \check{S}\,$ then 
$I:=(J\cap S)$ is a prime $\check{T}$ invariant ideal of $S$.   
\end{lem}
\begin{pf}
(i) The algebra $\check {S}$ admits a natural grading by $\check {T}$ through
$\check {S}_{t}:=S\otimes t$. Since $J$ is graded one can use ~\Lem{l2} (i).
Take homogeneous $a_1,a_2\in \check {S}\backslash J$. Write $a_i=s_it_i:\ 
s_i\in S,\ t_i\in T,\ i=1,2.$ Then $s_1,s_2\in S\backslash I$ so 
$s_2.t_1^{-1}\in S\backslash I$. Take
$g\in S$ such that $s_1g(s_2.t_1^{-1})\not\in I$. Then
$a_1(g.t_1)a_2=s_1g(s_2.t_1^{-1})t_1t_2\not\in J$ as required. 
The last part is clear.

(ii)
The adjoint action of $\check {T}$ defines a $H$ grading
on $\check {S}$ and on $S$. Since $\check {T}\subset\check{S}$ each two-sided 
ideal of $\check {S}$ is graded so $I$ is also graded. 
Assume that $I$ is not prime. Then, by~\Lem{l2},
there exist homogeneous $a,b\in S\setminus I$ such that
$aSb\subseteq I$. Then 
$\ a\check{S}b=aS\check{T}b=aSb\check{T}\subseteq I\check{T}\subseteq J\,$
that contradicts $J$ being prime and completes the proof.
\end{pf}
\section{some commutation relations in $R^w_0$}
\label{scr}
Fix $w\in W$. For a weight vector $a\in R^w$ denote by $\lwt a$ 
(resp., $\rwt a$) the weight of $a$ wrt the left (resp., right)
action of $T$.
If $L$ is a subspace of $R^w$ set $L\vert^{\lambda}=\{ a\in L: \lwt {a}=\lambda \}\, 
,\  L\vert_{\mu}=\{ a\in L: \rwt {a}=\mu \}.$ 
Given weight vector $\xi\in V^+(\lambda)\vert_{\mu}\ $ it is convenient to write
$\ c^{\lambda}_{\xi,u_{\lambda}}\,$ as $\ c^{\lambda}_{\mu}.$
\subsection{}
\label{scr1}
Let $J^+_{\lambda}(\eta)\,$ (resp., $J^-_{\lambda}(\eta)\,$) denote the left ideal
of $R^+$ generated by $c^{\lambda}_{\eta'}\,$ with $\, \eta' < \eta\,$ 
(resp., $\, \eta' > \eta\,$). In the notation of ~\cite{j}, 9.1.5 one has 
$\,J^{\pm}_{\lambda}(\eta)=J^{\pm}_{\lambda}(\eta,\lambda)\cap R^+ .$
By~\cite{j}, 9.1.5 $J^{\pm}_{\lambda}(\eta)\,$
are two-sided ideals of $R^+.$

The commutative relations ~\cite{j}, 9.1.5 imply that the following 
relations hold in $R^+\,$:

$$(i)\ \  c^{\nu}_{\mu} c^{\lambda}_{\eta}=q^{(\lambda , \nu )-(\eta , \mu)} 
c^{\lambda}_{\eta} c^{\nu}_{\mu}
  {\text {\ \ \ mod\ \ \ }} J^+_{\lambda}(\eta)|^{\lambda +\nu},$$   
$$(ii)\ \  c^{\lambda}_{\eta} c^{\nu}_{\mu}=q^{(\lambda , \nu )-(\eta , \mu)} 
c^{\nu}_{\mu} c^{\lambda}_{\eta}
  {\text {\ \ \ mod\ \ \ }} J^-_{\lambda}(\eta)|^{\lambda +\nu}.$$   

\subsection{}
\label{scr1w}
Let $J^+_{\lambda}(\eta)_w\ $ (resp., $J^+_{\lambda}(\eta)_w\ $)
denote the left ideal of $R^w_0$ generated 
by $c^{-\lambda}_w c^{\lambda}_{\eta'}\,$ with $\, \eta' < \eta \,$
(resp., $\, \eta' > \eta\,$).

\begin{lem}{jw}
For any $\ \lambda ,\nu \in P^+(\pi );\ \mu\in\Omega(V^+ (\nu )), 
\eta\in\Omega(V^+ (\lambda ))\ $ one has
$$(i)\ \ c^{-\lambda-\nu}_w c^{\nu}_{\mu} c^{\lambda}_{\eta}=
q^{(\lambda , \nu )-(\eta , \mu)} c^{-\lambda-\nu}_w  
c^{\lambda}_{\eta} c^{\nu}_{\mu}
{\text {\ \ \ mod\ \ \ }} J^+_{\lambda}(\eta)_w \, ,$$
$$(ii)\ \  c^{-\lambda-\nu}_w c^{\lambda}_{\eta} c^{\nu}_{\mu}=
q^{(\lambda , \nu )-(\eta , \mu)}  c^{-\lambda-\nu}_w 
c^{\nu}_{\mu} c^{\lambda}_{\eta}
{\text {\ \ \ mod\ \ \ }} J^-_{\lambda}(\eta)_w.$$
\end{lem}
\begin{pf} 
Consider $a\in J^+_{\lambda}(\eta)|^{\lambda +\nu}\, .$
By definition of $J^+_{\lambda}(\eta)$ one can write 
$a=\sum_i c^{\nu_i}_{\xi_i} c^{\lambda}_{\eta'_i}\ $, $\eta'_i<\eta$ for all $i.$
Since $\lwt a=\lambda +\nu\,$ one can assume that $\nu_i=\nu \,$ for all $i.$
Therefore
$$c^{-\lambda-\nu}_w a=\sum_i c^{-\lambda-\nu}_w c^{\nu}_{\xi_i} 
c^{\lambda}_{\eta'_i}=\sum_i b_i (c^{-\lambda}_w c^{\lambda}_{\eta'_i}),
 \text {\ \ \ where\ \ \ }  b_i\in R^w.$$
Since $\lwt {(c^{-\lambda-\nu}_w a)} =\lwt {(c^{-\lambda}_w c^{\lambda}_{\eta'_i})}=0$
it follows that $\lwt {b_i} =0$ for all $i$ so $b_i \in R^w_0.$ 

Consequently
$\ c^{-\lambda-\nu}_w J^+_{\lambda}(\eta)|^{\lambda +\nu}\subseteq 
J^+_{\lambda}(\eta)_w\ $ and similarly 
$\ c^{-\lambda-\nu}_w J^-_{\lambda}(\eta)|^{\lambda +\nu}\subseteq
J^-_{\lambda}(\eta)_w.$ 
Multiply relations (i), (ii) of ~\ref{scr1} on 
$c^{-\lambda-\nu}_w\, .$ Then the inclusions above
give the relations (i), (ii).
\end{pf}

\subsection{}
\label{autfi}
For $\nu\in P(\pi)\,$ consider the inner automorphism $\phi^{\nu}_w\,$ of $R^w$: 
$a\mapsto c_w^{-\nu} a c_w^{\nu}.$ Since $\phi^{\nu}_w\,$ preserves both
left and right weight subspaces its restriction on $R^w_0$ gives an automorphism
$\phi^{\nu}_w\,$ of $R^w_0$ which preserves the right weight subspaces.
Set $\Phi_w=\{\phi^{\nu}_w\ |\  \nu\in P(\pi)\}.$

From ~\cite{j}, 9.1.4(i), 10.1.11(ii) it follows that for weight vector $a\in R^w_0$
one has $\  a c_e^{\nu}=q^{(\nu, \rwt {a})} c_e^{\nu} a,\ 
a c_{w_0}^{\nu}=q^{(-w_0\nu, \rwt {a})} c_{w_0}^{\nu} a$. This implies that
$\, c_w^{-\nu}c_e^{\nu},\ c_w^{-\nu}c_{w_0}^{\nu}$ are normal elements of $R_0^w$
for all $\nu\in P^+(\pi)$. 

Take $\,\mu=w\nu.$ Then  ~\Lem{jw} gives
$$\text {(i)}\ \ \phi^{\nu}_w \left( c^{-\lambda}_w c^{\lambda}_{\eta}\right) =
q^{(w\nu, \eta -w\lambda)} c^{-\lambda}_w c^{\lambda}_{\eta}
{\text {\ \ \ mod\ \ \ }}  J^+_{\lambda}(\eta)_w,$$
Moreover $\  J^+_{\nu}(\mu )_w\ $ is $\Phi_w$-invariant.
$$\text {(ii)}\ \ \phi^{\nu}_w \left( c^{-\lambda}_w c^{\lambda}_{\eta}\right) =
q^{-(w\nu, \eta -w\lambda)} c^{-\lambda}_w c^{\lambda}_{\eta}
{\text {\ \ \ mod\ \ \ }}  J^-_{\lambda}(\eta)_w.$$
Moreover $\  J^-_{\nu}(\mu )_w\ $ is $\Phi_w$-invariant.

Let us show that the $\ J^{\pm}_{\lambda}(\eta)_w\  $ are two-sided ideals.
Take  $\ (c^{-\lambda}_w c^{\lambda}_{\eta'})\ $ with $\ \eta'<\eta.$ 
As noted in the proof of ~\Lem{jw} one has
$\ c^{-\lambda-\nu}_w J^+_{\lambda}(\eta)|^{\lambda +\nu}\subseteq 
J^+_{\lambda}(\eta)_w\ $ for any $\nu\in P^+(\pi).$ Then
$\ c^{-\lambda-\mu}_w c^{\lambda}_{\eta'}c^{\nu}_{\mu} \in 
J^+_{\lambda}(\eta)_w.$ Therefore  
$\ c^{-\lambda}_w c^{\lambda}_{\eta'} c^{\nu}_{\mu}c^{-\nu}_w=\phi^{-\nu}_w
\left( c^{-\lambda-\mu}_w c^{\lambda}_{\eta'}c^{\nu}_{\mu}\right)
\in J^+_{\lambda}(\eta)_w.$ Since the elements $\ c^{\nu}_{\mu}c^{-\nu}_w\ $
generate $R^w_0$ it follows that
$\ J^+_{\lambda}(\eta)_w\ $ is a two-sided ideal of $R^w_0$. 
The same reasoning applies to $\,J^-_{\lambda}(\eta)_w\, $.

Since the $\ J^{\pm}_{\lambda}(\eta)_w\  $ are two-sided $\Phi_w$-invariant ideals
and the $\ c_w^{\lambda}:\ \lambda\in P(\pi)$ generate $R^w$ over $R^w_0\,$
it follows that $\ R^w J^{\pm}_{\lambda}(\eta)_w=J^{\pm}_{\lambda}(\eta)_w R^w\ $
and $\ R^w J^{\pm}_{\lambda}(\eta)_w R^w\cap R^w_0=J^{\pm}_{\lambda}(\eta)_w.$

\subsection{}
\label{scr2}
\begin{lem}{}
For any $\ \lambda ,\nu \in P^+(\pi );\ \mu\in\Omega(V^+ (\nu )), 
\eta\in\Omega(V^+ (\lambda ))\ $ one has
$$\text {(i)}\ \ (c^{-\nu}_w c^{\nu}_{\mu}) (c^{-\lambda}_wc^{\lambda}_{\eta})=
q^{(\lambda , \nu )-(\mu ,w\lambda )} c^{-\lambda-\nu}_w c^{\nu}_{\mu}
c^{\lambda}_{\eta} \text {\ \ \ mod\ \ \ } J^+_{\nu}(\mu )_w,$$
$$\text {(ii)}\ \  (c^{-\nu}_w c^{\nu}_{\mu})(c^{-\lambda}_w c^{\lambda}_{\eta})=
q^{(\mu,\eta -w\lambda)} \phi^{\nu}_w \left( c^{-\lambda}_w c^{\lambda}_{\eta}\right)
(c^{-\nu}_w c^{\nu}_{\mu}) \text {\ \ \ mod\ \ \ } J^+_{\nu}(\mu )_w,$$
$$\text {(iii)}\ \ (c^{-\nu}_w c^{\nu}_{\mu}) (c^{-\lambda}_wc^{\lambda}_{\eta})=
q^{-(\lambda , \nu )+(\mu ,w\lambda )} c^{-\lambda-\nu}_w c^{\nu}_{\mu}
c^{\lambda}_{\eta} \text {\ \ \ mod\ \ \ } J^-_{\nu}(\mu )_w,$$
$$\text {(iv)}\ \   (c^{-\nu}_w c^{\nu}_{\mu})(c^{-\lambda}_w c^{\lambda}_{\eta})=
q^{-(\mu,\eta -w\lambda)} \phi^{\nu}_w \left( c^{-\lambda}_w c^{\lambda}_{\eta}\right)
(c^{-\nu}_w c^{\nu}_{\mu}) \text {\ \ \ mod\ \ \ } J^-_{\nu}(\mu )_w.$$
\end{lem}
\begin{pf}
(i) By~\ref{autfi} $\ c^{-\lambda}_wJ^+_{\nu}(\mu )_wc^{\lambda}_{\eta}\subseteq 
J^+_{\nu}(\mu )_w\, .$ Therefore, by~\ref{autfi}(i), one has
$$(c^{-\nu}_w c^{\nu}_{\mu}) c^{-\lambda}_wc^{\lambda}_{\eta}=c^{-\lambda}_w 
\phi^{-\lambda}_w (c^{-\nu}_w c^{\nu}_{\mu})c^{\lambda}_{\eta}=
q^{(\lambda , \nu )-(\mu ,w\lambda )} c^{-\lambda-\nu}_w c^{\nu}_{\mu}
c^{\lambda}_{\eta} \text {\ \  mod\ } J^+_{\nu}(\mu )_w .$$
The proof of (iii) is similar.

By~\Lem{jw}(i) one has 
$\ \ \ c_w^{-\lambda -\nu} c^{\nu}_{\mu}  
c^{\lambda}_{\eta}=q^{-(\nu ,\lambda)+(\mu ,\eta)}
c_w^{-\lambda -\nu} c^{\lambda}_{\eta }c^{\nu}_{\mu} \text { \ \ mod\ \ } 
J^+_{\nu}(\mu )_w\, .$

Taking into account the relation above the formula (i) takes the form
\begin{eqnarray*}
(c^{-\nu}_w c^{\nu}_{\mu})(c^{-\lambda}_w c^{\lambda}_{\eta})=
q^{(\lambda , \nu )-(\mu ,w\lambda )} c^{-\lambda-\nu}_w
c^{\nu}_{\mu} c^{\lambda}_{\eta}=q^{(\mu ,\eta)-(w\lambda ,\mu)}
c^{-\lambda-\nu}_w c^{\lambda}_{\eta} c^{\nu}_{\mu}=\\
q^{(\mu,\eta -w\lambda)} \phi^{\nu}_w \left( c^{-\lambda}_w c^{\lambda}_{\eta}\right)
(c^{-\nu}_w c^{\nu}_{\mu})\text {\ \ \ mod\ \ \ } J^+_{\nu}(\mu )_w .
\end{eqnarray*}
The proof of (ii) is similar.
\end{pf}
\section{spectral decomposition of $R^+,\ \check {R}^w_0,\ R^w_0$.}
\subsection{}
\label{fiinv}
The following construction is similar to~\cite{j}, 9.3.8.

Fix $\, P\in \Spec R^w_0$ or $\, P\in \Spec \check{R}^w_0$. 
For each $\nu\in P^+(\pi)$ set
$$C_P(\nu ):=\{ \mu\in \Omega (V(\nu ))\  |\ \exists \xi\in V(\nu )^*|_{\mu} :\ 
\  (c^{-\nu}_w c^{\nu}_{\xi})\not\in P \}.$$  
Obviously $w\nu\in C_P(\nu ).$ Denote by $D^+_P(\nu )$ (resp.,  $D^-_P(\nu )$)
the set of  minimal (resp., maximal) elements of $C_P(\nu ).$ 

Fix $\mu\in D^+_P(\nu ),\ a=(c^{-\nu}_w c^{\nu}_{\mu})\not\in P.$
Then $\, J^+_{\nu}(\mu )_w\subseteq P\,$ so, by ~\ref{scr2}(ii), one has
$$a (c^{-\lambda}_w c^{\lambda}_{\eta})= q^{(\mu,\eta -w\lambda)} 
\phi^{\nu}_w (c^{-\lambda}_w c^{\lambda}_{\eta}) a \text {\ \ \ mod\ \ \ } P.$$
Thus for homogeneous $\, b\in R^w_0\,$ one has
$$a b=q^{(\mu,\rwt b)}\phi^{\nu}_w (b) a \text {\ \ \ mod\ \ \ } P.$$

Thus $a$ is a normal element modulo $P$ and hence a non-zero divisor. 
It follows that if $b\ $ is homogeneous and $\, b\in P\,$ then 
$\phi^{\nu}_w (b)\in P.$ Thus we have proved the

\begin{lem}{crl}
Any $\check T$ invariant prime ideal of $R_0^w\,$ is 
$\Phi_w$ invariant.
\end{lem}

\subsection{}
\label{specplus}
Let $P^{++}\ $ be a set of regular dominant weights. Set 
$$R^{++}:=\sum _{\nu\in P^{++}} V^+(\nu),$$
$$\Spec_+ R^+:=\{P\in \Spec R^+\, :\  R^{++}\not\subseteq P\}.$$
In this subsection we will define a decomposition of $\Spec_+ R^+$.
\subsubsection{}
\label{dcm11}
Fix $\, P\in \Spec_+ R^+.$
Similar to ~\ref{fiinv} for each $\nu\in P^+(\pi)\ $ set
$$C_P(\nu ):=\{ \mu\in \Omega (V(\nu )) |\ \exists \xi\in V(\nu )^*|_{\mu}: 
\ c^{\nu}_{\xi}\not\in P \}.$$ 
Since $\ R^{++}\not\subseteq P\ $ it follows that
$C_P(\nu )\not =\emptyset\,$ for all  $\nu\in P^+(\pi).$ 
Denote by $D^+_P(\nu )$ (resp.,  $D^-_P(\nu )$)
the set of minimal (resp., maximal) elements of $C_P(\nu ).$ 
The reasoning in ~\cite{j}, 9.3.8 shows that
there exists $\ y_{\pm}\in W$ such that 
$D^{\pm}_P(\nu )=\{ y_{\pm}\nu \}.$ Denote by $\, X(y_-,y_+)\,$
the set of all $P\in \Spec_+ R^+\ $ such that
$D^-_P(\nu )=\{ y_-\nu \},\ \ D^+_P(\nu )=\{ y_+\nu \}.$ 
Since any $P\in X(y_-,y_+)$ contains $J^{\pm}_{\nu}(y_{\pm}\nu)$
for all $\nu\in P^+(\pi)$, the relations~\ref{scr1} imply 
that $c^{\nu}_{y_-}, c^{\nu}_{y_+}$ are normal modulo $P$.
\subsubsection{}
\label{tspect}
\begin{lem}{graded}
Take $P\in \Spec_+ R^+$. Then for all $\mu\in P(\pi)$ 
a subspace $P\cap R^+|^{\mu}$ (resp., $P\cap R^+|_{\mu}$) is 
graded wrt the right (resp., left) action of $T$.
\end{lem}
\begin{pf}
It is sufficient to check that for all $a\in (P\cap R^+|^{\mu})$
(resp., $a\in (P\cap R^+|_{\mu})$) one has $a.T\subset P$
(resp., $T.a\subset P$).
Take $y\in W$ such that $D^+_P(\nu )=\{ y\nu \}.$
Since $c^{\nu}_y$ is normal modulo $P$ we conclude from~\ref{scr1}(i)
that for any weight vector $c^{\lambda}_{\eta}$ and any $\nu\in P(\pi)$ 
one has
$$c^{\lambda}_{\eta}=q^{(\lambda,\nu)-(\eta,y\nu)}c^{\lambda}_{\eta}=
\tau(\nu).c^{\lambda}_{\eta}.\tau(y\nu) \text { mod } P.$$ 
Hence $a=\tau(\nu).a.\tau(y\nu)$ modulo $P$ for all $a\in R^+,\ \nu\in P(\pi)$.
If $a\in (P\cap R^+|^{\mu})$ then $\tau(\nu).a=q^{(\mu,\nu)}a$ so
$a.\tau(y\nu)\in P$. Similarly if
$a\in (P\cap R^+|_{\mu})$ then $\tau(\nu).a\in P$. 
This implies the required assertion. 
\end{pf}
\begin{rem}{ttspec}
The Lemma implies that the set of prime ideals of $R^+$ which are invariant wrt 
the left action of $T$ coincides with the set of primes
which are invariant wrt the right action of $T$. 
Therefore the same assertion holds for the ring $R^w$.
We will denote the corresponding sets of invariant ideals
by $(\Spec_+ R^+)^T,\ (\Spec R^w)^T$. 
\end{rem}
\subsubsection{}
\label{demaz}
Fix $y\in W.$ Denote by $V_y^{\pm}(\lambda)^{\perp}\,$ the
orthogonal of the Demazure module 
$\ V_y^{\pm}(\lambda):=U_q(\frak b ^{\pm}) u_{y\lambda}\,$ 
in $V(\lambda)^*\,$, the latter identified with $V^+(\lambda)\,$. 
Set
$$Q(y)^{\pm}:=\sum_{\lambda\in P^+(\pi)} V_y^{\pm}(\lambda)^{\perp}.$$
Observe that $Q(y)^{\pm}\supseteq J^{\pm}_{\nu}(y_{\pm}\nu)$
for all $\nu\in P^+(\pi)$ so $c^{\nu}_y$ is normal modulo
$Q(y)^{\pm}$. Observe also that $c_w\cap Q(y)^+=\emptyset$ 
(resp., $c_w\cap Q(y)^-=\emptyset$) if $w\leq y$ (resp., $w\geq y$).
 
By~\cite{j}, 10.1.8 $Q(y)^+\,$ is a completely prime ideal of $R^+\,$
(but note a slight difference of notation). 
A similar assertion holds for  $Q(y)^-.$ The reasoning in~\cite{j}, 10.1.13
shows that

\begin{prop}{dmz}
Every $P\in X(y_1,y_2)\,$ contains $Q(y_1)^-,\ Q(y_2)^+.$
\end{prop}

In particular, $\, Q(y_2)^+\,$ (resp., $Q(y_1)^-\, $) is a unique
minimal element of $X(e,y_2)\,$ (resp., $X(y_1,w_0)\,$).
\subsubsection{}
\label{dcm1}
The following lemma is a particular case of~\cite{j3}, 5

\begin{lem}{geom}
Let $P\in X(y_1,y_2),\ c^{\lambda}_y\not\in P\ $ for
some $\lambda\in P^{++}\, , y\in W.$ Then $\, y_1\leq y\leq y_2.$
\end{lem}
\begin{pf}
By~\Prop{dmz} $\ Q(y_2)^+\subseteq P\,$ so $\, c^{\lambda}_y\not\in Q(y_2)^+.$
The definition of  $Q(y_2)^+\, $ implies that
$\ u_{y\lambda}\in V_{y_2}^+(\lambda)\,$ so
$\ V_{y}^+(\lambda)\subseteq V_{y_2}^+(\lambda).$ By~\cite{j}, 4.4.5
it follows that $\, y\leq y_2.$ Similarly $\, y_1\leq y.$
\end{pf}

In particular, by the definition of $\, X(y_1,y_2)\,$, if $\ P\in X(y_1,y_2)\,$ 
then $\, c^{\lambda}_{y_1}\not\in P\ $. Therefore $\  y_1\leq y_2.$ Set
$$W\diamond W:=\{ (y_1,y_2)\in W\times W |\ y_1\leq y_2\}.$$

\begin{cor}{dcmrplus}
$$\Spec_+ R^+ =\coprod_{ (y_1,y_2)\in W\diamond W} X(y_1,y_2).$$
\end{cor}

\begin{rem}{}
It will be shown that each $X(y_1,y_2)$ is non-empty.
\end{rem}
\subsection{}
In this subsection we will define decompositions of $\Spec \check {R}^w_0,\ 
\Spec R^w_0$ which are similar to the above decomposition of $\Spec_+ R^+$.
\subsubsection{}
\label{cnn}
In order to relate $\,\Spec_+ R^+\,$
and $\,(\Spec R_0^w)^{\check {T}}$ recall that we have embeddings
\begin{equation}
R^+\overset{l_w}{\hookrightarrow}R^w\overset{\rho_0}{\hookleftarrow}R^w_0
\end{equation}
where  $\rho_0\,$ is the obvious embedding
and $l_w\,$ is the localization map.
For a two-sided ideal $I$ of $R^+\,$ (resp., of $R^w_0$) 
denote the ideal $R^w l_w(I)R^w\,$ (resp., $R^w \rho_0(I)R^w\,$) of
the ring $R^w\,$ by $I^l$ (resp.,  by $I^{\rho}$).

Let us show that the correspondence $I\mapsto I^{\rho}$ defines an 
order preserving injective map 
$$\rho:(\Spec R_0^w)^{\check{T}}\to(\Spec R^w)^{T}.$$

In fact, the torus $\{ c_w^{\nu }\}_{\nu\in P(\pi)}\ $ acts on $R_0^w$ by
automorphisms $\{\phi^{\nu}_w\}\ $ and  
$R^w=R_0^w\#\{ c_w\}$.
Let $P$ be a $\check {T}$ invariant prime ideal of $R_0^w$. Then, 
by~\Lem{crl},  $P$ is $\Phi_w$ invariant.
Then  $\ P^{\rho}=(P\# \{c_w\})\ $ is prime by~\Lem{l3}(i) and
is obviously $T$ invariant.
Moreover, $\, P^{\rho}\cap R_0^w=P.$ This gives an order preserving 
injection of
$\ (\Spec R_0^w)^{\check {T}}\,$ into $\, (\Spec R^w)^{T}.$

Furthermore, by~\cite{j}, A.2.8 and the noetherianity of $R^+$ (\ref{noeth}),
$l_w$ induces an  order preserving  bijection $P\mapsto P^l\,$ 
(with inverse $Q\mapsto Q\cap  R^+\,$) of
$\ \Spec_w R^+ :=\{ P\in \Spec R^+ |\ P\cap c_w=\emptyset\}\ $ onto $\Spec R^w.$
Since this bijection maps $T$ invariant prime ideals
to  $T$ invariant prime ideals, it induces
an order preserving injection of $\ (\Spec R_0^w)^{\check {T}}\ $ 
into $\ (\Spec_w R^+)^T$.
We may summarize the above by the following diagram:
\begin{equation}
(\Spec_w R^+)^{T}\overset{\sim}{\rightarrow} (\Spec R^w)^{T}
\overset{\rho }{\hookleftarrow}(\Spec R^w_0)
\label{diag1}
\end{equation}

\begin{rem}{checkq}
Let $Q\in \Spec_w R^+$ be a $T$ invariant completely prime ideal.
Then 
$$Q_w:=Q^l\cap R_0^w=\sum_{\lambda\in P^+(\pi)}c^{-\lambda}_w(Q\cap V^+(\lambda)).$$
is a  $\check T$ invariant completely prime ideal of 
$R^w_0$ so, by~\Lem{l3}(i), 
$\, \check {Q}_w:=(Q_w\#\check {T})\, $ is a 
completely prime ideal of $\check {R}_0^w.$
\end{rem}
\subsubsection{}
\label{tcomponents}
Fix $\ P\in(\Spec R_0^w)^{\check T}\ $  and set $\, P'=(P^{\rho}\cap R^+).$

Since $\ P^{\rho}\cap R^w_0=P\ $ it follows that 
$(c^{-\nu}_w c^{\nu}_{\xi})\in P\,$ iff $\, c^{\nu}_{\xi}\in P'.$
Therefore $\, J^{\pm}_{\nu}(\mu )_w\subseteq P\,$
iff $\, J^{\pm}_{\nu}(\mu )\subseteq P'.$
Hence $\,D^{\pm}_P(\nu)=D^{\pm}_{P'}(\nu)\,$ for all $\,\nu\in P^+(\pi).$ 

Since $\, P'\in \Spec_w R^+\subset\Spec_+ R^+ $ there exist $\, y_{\pm}\in W\,$ 
such that $D^{\pm}_P(\nu)=D^{\pm}_{P'}(\nu )=\{ y_{\pm}\nu \}.$ 
Since $\,P'\cap c_w=\emptyset$,
we conclude from~\Lem{geom} that $\, y_-\leq w\leq y_+.$
\subsubsection{}
\label{components}
Fix $P\in\Spec R_0^w$ (resp., $P\in\Spec \check{R}_0^w$) and let $P'$ be a 
maximal $\check T$ invariant ideal contained in $P$ (resp., $P'=P\cap R^w_0$). 
Then $\,D^{\pm}_P(\nu)=D^{\pm}_{P'}(\nu)\,$ for all $\,\nu\in P^+(\pi).$ 
By~\Lem{l3} $P'\in(\Spec R_0^w)^{\check T}$. Hence
$D^{\pm}_P(\nu)=\{ y_{\pm}\nu \}$ for some $y_{\pm}$ such that
$\, y_-\leq w\leq y_+.$ Set
$$W\overset{w}{\diamond}W:=\{ (y_1,y_2)|\ y_1\leq w\leq y_2 \}.$$

Fix $\ (y_1,y_2)\in W\overset{w}{\diamond}W\ $ and let $X_w(y_1,y_2)$ 
(resp., $Y_w(y_1,y_2)$) denote the set of all $P\in \Spec \check{R}_0^w$
(resp., $P\in\Spec R_0^w$) such that
$\ D^-_P(\nu )=\{ y_1\nu \},\ D^+_P(\nu )=\{ y_2\nu\}\ $ for all 
$\ \nu\in P^+(\pi).$ 
We summarize the results above by the

\begin{prop}{dcms}
$$\text{(i)  }\Spec \check{R}_0^w=\coprod_{(y_1,y_2)\in W\overset{w}{\diamond}W} 
X_w(y_1,y_2).$$
$$\text{(ii)  } \Spec R_0^w=\coprod_{(y_1,y_2)\in W\overset{w}{\diamond}W} 
Y_w(y_1,y_2).$$
\end{prop}

\section{the study of the strata}
The goal of this section is to show that for each $\, (y_1,y_2)\in W
\diamond W\,$ the component  $X(y_1,y_2)\,$  of $\,\Spec R^+\,$ 
has a unique minimal element $\, Q(y_1,y_2).$ Moreover for 
$y_1\leq w\leq y_2$ the ideals $\,Q(y_1,y_2)_w,\ \check {Q}(y_1,y_2)_w$
(notations of ~\Rem{checkq}) are unique minimals of $Y_w(y_1,y_2),X_w(y_1,y_2)$
respectively.

\subsection{Notations}
\subsubsection{}
\label{epsilon}
Set $\, U:=U_q(\frak g).$ For $i=1,\ldots ,l\ $ set 
$\varphi _i (a):=\max\{n: a.y_i^n\not =0\}$ 
(resp., \\ $\varepsilon _i (a):=\max\{n: a.x_i^n\not =0\}$)
for all $a\in R^+$ non-zero; also set $\varphi _i (0):=0,\
\varepsilon _i (0):=0.$ Note that
$$ \varphi _i (ab)= \varphi _i (a) +\varphi _i (b)\ \ \text {for non-zero a,b,\ \ }$$
$$(\alpha _i,\rwt a)=\varphi _i (a)-\varepsilon _i (a)\text
 {\ \ for any weight vector \ }a.$$

Let $a\in R^+\,$ be a non-zero weight vector.
Define $\ a.y_i^*:=a.y_i^{\varphi _i (a)}\,$ 
(resp., $\ a.x_i^*:=a.x_i^{\varepsilon _i (a)}$).
Furthermore for a fixed reduced decomposition $\, w=s_{i_1}\ldots s_{i_r}\,$
(resp., $w w_0 =s_{j_1}\ldots s_{j_p}\ $)
set $\, a.y_w^*:=a.y_{i_1}^*\ldots y_{i_r}^*\, $
(resp., $\, a.x_w^*:=a.x_{j_1}^*\ldots x_{j_p}^*$).

Recall that $V^+(\nu)\cong V(\nu)^*$ as right $U$ modules for all 
$\nu\in P^+(\pi)$.
In particular $V^+(\nu)\,$ has highest weight $\nu$ and the corresponding
highest weight vector is annihilated by the $y_i:\ i=1,\ldots ,l$
rather than by the $x_i$. Moreover $\ \varepsilon _i (c_w^{\nu})=0$ (resp., 
$\,\varphi _i(c_w^{\nu})=0$) if $s_iw<w\, $ (resp., if $s_iw>w$).
It implies that $\, c^{\nu}_w.y_w^*=c^{\nu}_e,\ \ 
c^{\nu}_w.x_w^*=c^{\nu}_{w_0}\,$ up to non-zero scalars.

Fix $\, i\in \{ 1,\ldots ,l\}.$
Suppose $a,b$ are weight vectors and set $\ \varphi _i (a)=n,\ \varepsilon_i (a)=n',\   
\varphi _i (b)=m,\  \varepsilon_i (b)=m'.$ Since
$$\triangle (y_i)=y_i\otimes 1+ t_i\otimes y_i,\ \ \  
\triangle (x_i)=x_i\otimes t_i^{-1}+1\otimes x_i$$ 
it follows that there exist $\ P_{m+n}^n\in K^*\ $ such that $P_{m+n}^n=P_{m+n}^m\,$
and
$$(ab).y_i^*=P_{m+n}^n q^{(m\alpha_i ,\rwt a )} (a.y_i^*)( b.y_i^*)\, ,\ \ \ \ 
(ab).x_i^*=P_{m'+n'}^{n'} q^{-(n'\alpha_i ,\rwt b )} (a.x_i^* )(b.x_i^* ).$$

\subsubsection{}
\label{ntt}
Fix $w\in W$. Using notations of~\ref{demaz} set 
$$Q(y)^{\pm}_w:=\sum_{\nu\in P^+(\pi)} c_w^{-\nu} V_y^{\pm}(\nu)^{\perp}.$$
The ideal $Q(y)_w^+$ (resp., $Q(y)_w^-$) does not coincide with whole $R^w_0$
iff $y\geq w$ (resp., $y\leq w$); in this case, by~ \Rem{checkq}, 
it is a $\check T$ invariant completely prime ideal of $R^w_0$.

Recall that $\phi^{\nu}_w:\ a\mapsto c_w^{-\nu} a c_w^{\nu}\ $ is an 
automorphism of $ R^w\,$ and of $R^w_0.$ By~\Lem{crl} 
$\, Q(y)^{\pm}_w\ $ are $\Phi_w$ invariant. 
\subsubsection{}
\begin{defn}{sf} 
Fix $w\in W$. For $\eta\in wQ^- (\pi)\,$ call $\lambda\in P^+(\pi)\,$
sufficiently large for $\eta\,$  if the natural embedding
$\ c_w^{-\lambda} V(\lambda)^+|_{w\lambda +\eta}\hookrightarrow R^w_0|_{\eta}\ $
is bijective.
\end{defn}
Since $\dim R^w_0|_{\eta}<\infty\, $ the existence of such $\lambda\ $
follows from ~(\ref{A}).

\subsection{}
\label{identification}
\begin{lem}{} 
Take $\ \eta\in wQ^- (\pi)\,$ and choose $ \lambda\,$ sufficiently large for $\eta.$ 
Then $\, V^+(\lambda )|_{w\lambda +\eta}\ $ is $\,\Phi_w$ invariant.
\end{lem}
\begin{pf}
Identify the vector spaces  $\ R^w_0|_{\eta}$ and
$V^+(\lambda )|_{w\lambda +\eta}\ $
through the map $\, a\mapsto c_w^{\lambda} a.$

An automorphism $\,\phi_w^{\nu}\,$ leaves $\ R^w_0|_{\eta}\,$ invariant.
Then for any $\, a\in  R^w_0|_{\eta}$ one has 
$$\phi_w^{\nu}(c_w^{\lambda}a)=c_w^{-\nu}(c_w^{\lambda}a)c_w^{\nu}=
c_w^{\lambda}\phi_w^{\nu}(a)\in c_w^{\lambda} R^w_0|_{\eta}=
V^+(\lambda )|_{w\lambda +\eta}.$$
\end{pf}
\begin{rem}{identify}
Actually we showed that the bijection between  $\,R^w_0|_{\eta}\, $
and $\, V^+(\lambda )|_{w\lambda +\eta}\ $ commutes with the action
of $\Phi_w.$
\end{rem}
\subsection{}
Fix $\ \eta\in wQ^- (\pi)\,$ and choose $ \lambda\,$ sufficiently large for 
$\eta.$ 
Let us show that the eigenvalues of $\phi^{\nu}_w\,$ on $R^w_0|_{\eta}\ $ are  
some integer powers of $q.$ For this we will identify
$R^w_0|_{\eta}\ $ with  $\, V^+(\lambda )|_{w\lambda +\eta}\ $
and will study the change of the eigenvalues when we pass from
$\,\phi^{\nu}_w\,$ to $\,\phi^{\nu}_{s_i w}.$ 

Let $\overline {K}$ be the algebraic closure of $K.$ 
Set $\overline {V}^+(\lambda)=V^+(\lambda)\otimes_K\overline {K}.$
\subsubsection{}
\begin{lem}{intvl} Fix $\ \nu, \lambda\in P^+(\pi)$. Suppose  
$c_{\xi}^{\lambda}\in \overline {V}^+(\lambda )\,$ is a weight
vector such that
$$(a)\ \ \ (\phi^{\nu}_w)^m(c_{\xi}^{\lambda})\in \overline {V}^+(\lambda )\ 
\text {\ \ \ for all \ } m \in\Bbb {N},$$ 
$$(b)\ \ \ (\phi^{\nu}_w-s\cdot\id)^r(c_{\xi}^{\lambda})=0 \text {\ \  for some\ }
s\in\overline {K},\ r\in\Bbb {N}.$$ Then

(i) If $i\in \{ 1,\ldots ,l\}$ is such that $s_i w<w\,$ then
$$\ (\phi^{\nu}_{s_iw})^m(c_{\xi}^{\lambda}.y_i^*)\in 
\overline {V}^+(\lambda )\text 
{\ \ for all \ } m \in\Bbb {N} \text {\ \ and }$$
$$(\phi^{\nu}_{s_iw}-s'\cdot\id)^r(c_{\xi}^{\lambda}.y_i^*)=0 \text {\ \ where\ }
s'=s\cdot q^{(\rwt \xi ,w\nu)-(\rwt {(\xi .y^*_i)},s_i w\nu)}.$$

(ii) If $i\in \{ 1,\ldots ,l\}$ is such that $s_i w>w\,$ then
$$\ (\phi^{\nu}_{s_iw})^m(c_{\xi}^{\lambda}.x_i^*)\in \overline {V}^+
(\lambda )\text 
{\ \ for all \ } m \in\Bbb {N} \text {\ \ and }$$
$$(\phi^{\nu}_{s_iw}-s'\cdot\id)^r(c_{\xi}^{\lambda}.x_i^*)=0 \text {\ \ where\ }
s'=s\cdot q^{-(\rwt \xi ,w\nu)+(\rwt {(\xi .x^*_i)},s_i w\nu)}.$$
\end{lem}
\begin{pf}
We prove (i) by induction on the nilpotence degree $r$. 
Fix $i$ and set $\,\varphi:=\varphi _i,\ y:=y_i, \ m:=\varphi (c_w^{\nu}).$ 
Since $s_i w<w\,$ it follows from~\ref{epsilon} that 
$\, c_w^{\nu}.y^m=c_{s_i w}^{\nu}\,$ up to a non-zero scalar.

Set $c_{\xi _1}^{\lambda}:=(\phi^{\nu}_w-s\cdot\id)(c_{\xi}^{\lambda}).$ Then
$\,(\phi^{\nu}_w-s\cdot\id)^{r-1}(c_{\xi_1}^{\lambda})=0\ $ and also
$\ (\phi^{\nu}_w)^m(c_{\xi}^{\lambda})\in \overline {V}^+(\lambda )\ 
\text {\ for all \ }\\ m \in\Bbb {N}.$
One has
$\phi_w^{\nu} (c_{\xi}^{\lambda})=s c_{\xi}^{\lambda} +c_{\xi _1}^{\lambda}\ $
or, in other words,
\begin{equation}
c_{\xi}^{\lambda} c_w^{\nu}=s c_w^{\nu} c_{\xi}^{\lambda}+
c_w^{\nu} c_{\xi _1}^{\lambda}. 
\label{B}
\end{equation}
If $\, r=1\,$ then $\,\xi_1=0\,$ otherwise  $\rwt {\xi}=\rwt {\xi_1}.$

Set $n:=\varphi (c_{\xi}^{\lambda}),
\ n_1:=\varphi (c_{\xi_1}^{\lambda}).$ 
Then $\varphi  (c_{\xi}^{\lambda} c_w^{\nu})=m+n,\ 
\varphi  (c_{\xi_1}^{\lambda} c_w^{\nu})=m+n_1.$ From the formula ~(\ref{B})
it follows that $\ m+n_1\leq m+n. $ Therefore $n_1\leq n.$

Act by $y^{m+n}\,$ on the both sides of ~(\ref{B}). 
Applying ~\ref{epsilon} we get 
\begin{equation}
q^{(m\alpha ,\rwt \xi)}(c_{\xi}^{\lambda}.y^*) 
c_{s_i w}^{\nu}=q^{(n\alpha ,w\nu)}(s c_{s_i w}^{\nu} (c_{\xi}^{\lambda}.y^*)
+c_{s_i w}^{\nu} (c_{\xi_1}^{\lambda}.y^n)).
\label{C}
\end{equation}
Note that
$$(\rwt \xi ,w\nu)-(\rwt {(\xi .y^*_i)},s_i w\nu)=(\rwt \xi ,w\nu)-
(\rwt \xi +n\alpha,w\nu +m\alpha)=$$ 
$$ -(n\alpha ,s_i w\nu)-
(m\alpha ,\rwt \xi)=(n\alpha ,w\nu)-
(m\alpha ,\rwt \xi).$$
Therefore from the formula ~(\ref{C}) it follows that
\begin{equation}
(\phi^{\nu}_{s_iw}-s'\cdot\id)(c_{\xi}^{\lambda}.y^*)=(s'/s) 
c_{\xi_1}^{\lambda}.y^n.
\label{D}
\end{equation}

Since $\xi_1=0$ for $r=1$, the assertion for
this case immediately follows from ~(\ref{D}).

Suppose $\ n_1<n.$ Then $\ c_{\xi_1}^{\lambda}.y^n=0\,$ so the assertion 
holds.
Finally, if $\ n_1=n\,$ then $\ c_{\xi_1}^{\lambda}.y^n=
\ c_{\xi_1}^{\lambda}.y^*\ $
and $\,\rwt {(\xi_1.y^*)}=\rwt {(\xi.y^*)}.$
The induction hypothesis implies that
$$(\phi^{\nu}_{s_iw}-s'\cdot\id)^{r-1}(c_{\xi_1}^{\lambda}.y^*)=0,\ \ \ 
(\phi^{\nu}_{s_iw})^m(c_{\xi_1}^{\lambda}.y^*)\in \overline {V}^+(\lambda )\ 
\text {\ for all \ } m \in\Bbb {N}.$$ 
taking into account ~(\ref{D}) we get the required assertion.
The proof of (ii) is completely similar.
\end{pf}
\subsubsection{}
\label{intval}
By~\cite{j}, 9.1.4(i), 10.1.11(ii) one has  
$$\ c^{-\nu}_e c^{\lambda}_{\mu} c^{\nu}_e=q^{(\nu ,\mu- \lambda )} 
c^{\lambda}_{\mu} ,\ \ \ \  
\ c^{-\nu}_{w_0} c^{\lambda}_{\mu}c^{\nu}_{w_0}=
q^{-(w_0\nu ,\mu-w_0\lambda ))} c^{\lambda}_{\mu} .$$  

So all eigenvalues of the automorphisms $\phi_e^{\nu},\ \phi_{w_0}^{\nu}\ $ 
are integer powers of $q.$ Then from~\Lem{intvl} it follows, by induction, 
that for any $w\in W\,$ all eigenvalues of the automorphisms
$\phi_w^{\nu}\ $ are integer powers of $q.$ 

\subsection{}
\label{defwtw}
Since all eigenvalues of the system of automorphisms 
$\Phi_w \,$ are integer powers of $q$ it
follows that for each common eigenvector $a\in R^w$
there exists $\mu\in Q(\pi )\,$ such that
$\, \phi_w^{\nu}(a)=q^{(\mu , \nu)} a.$ This element $\mu\in Q(\pi )\,$
will be called eigenvalue of $\Phi_w$. From this we make the

\begin{defn}{wtw}
For $a\in R^w\ $ set $\wt_w a:=\mu\in Q(\pi )\  $ if
$\ \ \forall\  \nu \ \exists\  r\in \Bbb {N}:\  
(\phi_w^{\nu}-q^{(\mu , \nu)} \id )^r a=0.$
\end{defn}

\subsubsection{}
\label{gneigvl}
Suppose $a\in R^+$ is homogeneous and $\wt_w a\ $ is defined.
Then by~\Lem{intvl} $\wt_{s_iw} (a.y_i^*)$
(resp., $\wt_{s_iw} (a.x_i^*)$) is defined for $s_iw<w$ (resp., $s_iw>w$) 
and satisfies to the following relations:$$\left\{ 
\begin{array}{ll}
\wt_w a+w^{-1}\rwt a=\wt_{s_iw} (a.y_i^*)+(s_iw)^{-1}\rwt (a.y_i^*)
\text {\ \ \ if\ \ \ } s_iw<w \\
\wt_w a-w^{-1}\rwt a=\wt_{s_iw} (a.x_i^*)-(s_iw)^{-1}\rwt (a.x_i^*)
\text {\ \ \ if\ \ \ } s_iw>w 
\end{array}\right.$$
By induction for any reduced decomposition of $w\,$ (resp., $ww_0$) 
$$\wt_w a+w^{-1}\rwt a=\wt_e (a.y_w^*)+\rwt (a.y_w^*),
\ \ \ \ 
\wt_w a-w^{-1}\rwt a=\wt_{w_0} (a.x_w^*)-w_0\rwt (a.x_w^*).$$
The relations ~\ref{intval} imply that 
$$\wt_e c_{\xi }^{\lambda}=
\rwt {(c_e^{-\lambda} c_{\xi }^{\lambda})},\ \ \ \ 
\wt_{w_0} c_{\xi }^{\lambda}=
-w_0\rwt {(c_{w_0}^{-\lambda} c_{\xi }^{\lambda})}.$$

Hence one has the

\begin{prop}{peigval}
Take a weight vector $\ c_{\xi }^{\lambda}$ such that 
$\,\wt_w c_{\xi }^{\lambda}\,$ is defined. Then 
$$\wt_w c_{\xi }^{\lambda}+
w^{-1}\rwt {(c_w^{-\lambda} c_{\xi}^{\lambda})}
=2\rwt {(c_e^{-\lambda} c_{\xi .y_w^*}^{\lambda})} ,\ \ \ 
\wt_w c_{\xi }^{\lambda}-
w^{-1}\rwt {(c_w^{-\lambda} c_{\xi}^{\lambda})}
=-2w_0\rwt {(c_{w_0}^{-\lambda} c_{\xi .x_w^*}^{\lambda})} .$$
\end{prop}

Consider $\ a\in R^w_0|_{\eta}\ $ such that $\wt_w a\,$ is defined.
Note that $\wt_w a=\wt_w (c^{\lambda}_w a)\ $ for all $\,\lambda\in P(\pi).$ 
Choose $\lambda$ sufficiently large for $\eta$ (\Defn{sf}) and set  
$\ c_{\xi }^{\lambda}:=c^{\lambda}_w a.$
Then from the proposition above we get that
 
\begin{equation}
\left. 
\begin{array}{rr}
(\wt_w a+w^{-1}\eta)=2\rwt {(c_e^{-\lambda} c_{\xi .y_w^*}^{\lambda})}
\in 2 Q^-(\pi )\\ 
(\wt_w a-w^{-1}\eta)=-2w_0\rwt {(c_{w_0}^{-\lambda} 
c_{\xi .x_w^*}^{\lambda})} \in 2 Q^+(\pi )
\end{array}
\right\}\Longrightarrow
{w^{-1}\eta\leq\wt_w a\leq -w^{-1}\eta.}
\label{0}
\end{equation}
Note that $w^{-1}\eta\in Q^-(\pi )$.
\subsection{}
\label{twistaut}
Fix $w\in W$. Consider a twisted system of automorphisms 
$\tilde\Phi_w:=\{\tilde{ \phi}_w^{\nu}\}\,$ of $R^w_0\,$ given by 
$$ a\mapsto q^{(w^{-1}\rwt a,\nu)} \phi_w^{\nu} (a),\text {   
on any weight vector } a.$$
Since $J^+_{\nu}(w\nu)_w\subset Q(w)^+_w\ $ for any 
$\nu\in P^+(\pi )$, we conclude
from ~\Lem{jw}(i) that for any weight vector $a\in R^w_0$ one has
$\phi_w^{\nu}(a)=q^{(\nu,-w^{-1}\rwt a)}a \text {  mod  }Q(w)^+_w$. Therefore
\begin{equation}
\tilde{\phi}_w^{\nu}(a)=a\ \text {  mod  } Q(w)^+_w\ \ \text { for all } a\in R^w_0.
\label{QW}
\end{equation}
For each $\mu\in Q(\pi)$ denote by $L(w,\mu)|_{\eta}$ the maximal 
subspace of $R^w_0|_{\eta}$ on which all the endomorphisms
$\,(\tilde{ \phi}_w^{\nu}-q^{(\nu,\mu)}\id)\, ,\nu\in P(\pi)\,$ act
nilpotently. Set $L(w,\mu):=\underset{\eta}{\oplus} L(w,\mu)|_{\eta}.$
One has
\begin{equation}
L(w,\mu)=\sum \{ a\in R^w_0|\ 
\wt_w a=\mu-w^{-1}\rwt a\} .
\label{L}
\end{equation}
Then~(\ref{0}) implies that
$$R^w_0=\underset {\mu\in 2Q^-}{\oplus} L(w,\mu ).$$
Observe that $L(w,\mu)L(w,\nu)\subseteq L(w,\mu+\nu)\ $
so $L(w,0)$ is a subalgebra of $R^w_0$. Set 
$L'(w):=\underset{\mu\not= 0}{\oplus} L(w,\mu).$
\subsubsection{}
\begin{lem}{qww}
(i) One has $Q(w)^+_w=L'(w)$. In particular $R^w_0=L(w,0)\oplus Q(w)^+_w$.

(ii) Take a weight vector $c^{\lambda}_{\xi}$ such that
$\wt_wc^{\lambda}_{\xi}$ is defined. Then 
$$c^{\lambda}_{\xi}\in Q(w)^+\ \ \Longleftrightarrow\ \ 
\wt_w c^{\lambda}_{\xi}+w^{-1}\xi-\lambda\not=0.$$ 
\end{lem}
\begin{pf}
(i) Fix $\mu\not=0$ and $\nu\in P^+(\pi)$ such that $(\nu,\mu)\not=0$.
Take $a\in L(w,\mu)$. 
Since $(\tilde{\phi}_w^{\nu}-q^{(\nu,\mu)}\id)^r(a)=0$
for some $r\in\Bbb N$, we conclude from the 
formula ~(\ref{QW}) that $a\in Q(w)^+_w$. Hence $L'(w)\subseteq Q(w)_w^+$.

Now suppose that $ Q(w)_w^+\not\subseteq L'(w).$
The formula ~(\ref{QW}) implies that $Q(w)_w^+\, $ is $\tilde\Phi_w$ invariant.
Then there exists a weight vector
$a\in  Q(w)_w^+$ such that $a\in L(w,0).$ 
Since each automorphism $\tilde{ \phi}_w^{\omega_i}$ acts on
$L(w,0)$ nilpotently one can assume that $a$ is an
eigenvector that is $\ \tilde{ \phi}_w^{\nu}(a)=a\,$
for all $\nu\in P(\pi).$ Choose $ \lambda $ sufficiently
large for $\rwt a $ and write $\ a=c_w^{-\lambda} c_{\xi }^{\lambda}.$
 
From ~\Prop{peigval} and the definition
of  $\,\tilde{ \phi}_w^{\nu}\,$ we conclude
that $\ \rwt(c_e^{-\lambda}c_{\xi.y^*_w}^{\lambda})=0\,$ and so 
$c_{\xi.y^*_w}^{\lambda}=c^{\lambda}_{\lambda}\,$ up to a non-zero scalar.
Therefore
$$0\not=\xi.y^*_w(v_{\lambda})=\xi.(y_{i_1}^{n_1}\ldots y_{i_r}^{n_r})
(v_{\lambda})=\xi (y_{i_1}^{n_1}\ldots y_{i_r}^{n_r}v_{\lambda}).$$

By~\cite{j}, 4.4.6 $(y_{i_1}^{n_1}\ldots y_{i_r}^{n_r}v_{\lambda})\in
V_w^+(\lambda)\ $ so $\xi (V_w(\lambda)^+)\not= 0.$

However $\ a=c_w^{-\lambda} c_{\xi }^{\lambda}\in  Q(w)^+_w\ $ that is
$c_{\xi }^{\lambda}\in Q(w)^+.$ Hence
$\xi (V_w(\lambda)^+)=0\,$ giving the required contradiction.

(ii) Recall that $\, c^{\lambda}_{\xi}\in Q(w)^+\ $ iff 
$c^{-\lambda}_wc^{\lambda}_{\xi}\in Q(w)^+_w$. Then (i) 
and~(\ref{L}) imply
the required assertion.
\end{pf}
\begin{rem}{rssaction}
The lemma above and the formula ~(\ref{QW}) imply that $\tilde{\phi}_w^{\nu}(a)=a$ for 
all $\nu\in  P^+(\pi)$ iff $a\in L(w,0)$.
\end{rem}

\subsection{}
\label{qywy}
\begin{lem}{qyw}
$Q(y,w)_w:=Q(w)_w^+ +Q(y)_w^-\ $ is a completely prime ideal of $R_0^w$ for all  
$y\leq w$.
\end{lem}
\begin{pf}
By~\Lem{crl} $\, Q(y)_w^-\,$ is $\Phi_w$ invariant so $\tilde\Phi_w$ 
invariant. By~\Lem{qww}(i) $\ L'(w)=Q(w)_w^+\ $ therefore
$$Q(y,w)_w=L'(w)\oplus (L(w,0)\cap  Q(y)_w^-).$$
Consequently,
$$R^w_0/ Q(y,w)_w=(L(w,0)\oplus L')/ \left( (L(w,0)\cap Q(y)_w^- )\oplus
 L'\right)\cong L(w,0)/ (L(w,0)\cap Q(y)_w^- ).$$

To show that $\ L(w,0)/ (L(w,0)\cap Q(y)_w^- )\ $ is a domain,
observe that, by~\ref{ntt}, $\, Q(y)^-_w\ $ is a completely prime
ideal of $R^w_0.$ Since $L(w,0)$ is a subalgebra of $R^w_0$ it follows 
that $\ (L(w,0)\cap Q(y)_w^-)\ $ is a completely prime 
ideal of $L(w,0).$
\end{pf}
\subsubsection{}
\label{tilde}
Similar to ~\ref{twistaut} one can consider a twisted system of 
automorphisms $\{\tilde{ \phi}_w^{\nu}\}\ $
of $R^w_0$ given by $\ a\mapsto q^{-(w^{-1}\rwt a,\nu)} \phi_w^{\nu} 
(a),$ on any weight vector $a.$ Then reasoning similar 
to ~\ref{qww}---~\ref{qywy} shows that $Q(w,y)_w\,$ is
a completely prime ideal of the ring $R^w_0$ for all $y\geq w$.

\subsection{}
\label{qinr}
Fix $(y,w)\in W\diamond W$. By~\ref{demaz}
every $P\in X(y,w)$ contains $\tilde {Q}(y,w):=(Q(y)^{-}+Q(w)^{+})$.
The ideal $\tilde {Q}(y,w)$ is not in general prime. We describe now an operation
which, being applied to $\tilde {Q}(y,w)$, gives a prime ideal.

Recall that for all $z\in W$ the set $c_z$ is an Ore set in $R^+$. Let $I$ be a 
two-sided ideal in $R^+$ such that $I\cap c_z=\emptyset$. We define the 
saturation of $I$ along $c_z$ by the formula
$$I:c_z=\Ker\left( R^+\rightarrow(R^+/I)[c_z^{-1}]\right) .$$
For all $\nu\in P^+(\pi)$ the $c^{\nu}_w$ is normal modulo $\tilde {Q}(y,w)$ 
and modulo any $P\in X(y,w)$. Therefore $P:c_{w}=P$. Since the 
saturation along $c_w$ preserves the inclusion relation of ideals, it follows that 
$\ P\supseteq \tilde {Q}(y,w):c_w\ $ for all $P\in X(y,w)$. Set
$$Q(y,w):=\tilde {Q}(y,w):c_w=\{a\in R^+ |\ \exists\lambda\in P^+(\pi)\ 
s.t. \  c_w^{\lambda}a\in Q(y)^- +Q(w)^+\}.$$
Therefore $Q(y,w)=R^w Q(y,w)_w\cap R^+$.
By~\Lem{qyw} $Q(y,w)_w\, $ is a $\check {T}$ invariant completely prime 
ideal of $R^w_0.$ By~\ref{cnn} this implies that $Q(y,w)$ is a 
$T$ invariant completely prime ideal of $R^+.$
\subsection{}
\label{min}
\begin{prop}{mini}
The $T$ invariant  completely prime ideal $\,Q(y,w)\,$ of $R^+$
is the unique minimal element of $X(y,w)\,$ for all $\ (y, w)\in  W\diamond W.$ 
\end{prop}
\begin{pf}
By~\ref{qinr} any $P\in  X(y,w)\,$ contains $\,Q(y,w)$, which 
is a $T$ invariant completely prime ideal of $R^+.$ Therefore
it is sufficient to show that $\,Q(y,w)\in X(y,w).$

Recall that  
$$Q(y,w)=\{a\in R^+ |\ \exists\lambda\in P^+(\pi)\ s.t. \ 
 c_w^{\lambda}a\in Q(y)^- +Q(w)^+\}.$$
Since $\ c_w\cap  Q(y,w)=\emptyset\, , $ it suffices
to check that  $\  c_y^{\nu}\not\in Q(y,w)\ $
for all $\ \nu\in P^+(\pi).$ We prove this by induction. Namely,
from the pair $\ (y, w)\in  W\diamond W\ $ such that 
$c_y^{\nu}\in Q(y,w)\ $ we will construct a pair
$\ (s_iy, w')\in  W\diamond W\ $ such that $s_iy>y$ and
$\  c_{s_iy}^{\nu}\in Q(s_iy,w').$ Note that $(w_0,z)\in  W\diamond W\ $
forces $z=w_0$. Since $\  c_{w_0}^{\nu}\not\in Q(w_0,w_0)\ $
we will thus obtain a contradiction. The required assertion is 
proved in~\ref{aff}---~\ref{mini6} below.

\subsubsection{}
\label{aff}
Suppose that there exists $\nu\in P^+(\pi)\,$
such that $\,c_y^{\nu}\in Q(y,w).$
Then $c_w^{-\nu}c_y^{\nu}\in Q(y,w)_w$.
Set $\eta:=(y\nu - w\nu).$ By~\ref{twistaut} and the proof of~\Lem{qyw} 
one can write
\begin{equation}
c_w^{-\nu}c_y^{\nu}=\sum_{j=0}^m b_j,
\label{b}
\end{equation}
where the $b_j$ are weight vectors of the weight $\eta$, the values $\wt_w b_j$
are defined and pairwise distinct, $b_0\in L(w,0)\cap Q(y)^-_w\ $ 
and $b_j\in Q(w)^+_w\ $ for $j=1,\ldots,m$.

Choose $\mu \,$ sufficiently large for $\eta$ (see ~\Defn{sf}) such
that 
$\mu>\nu$ and set $\ \lambda:=\mu-\nu.$ For $i=0,\ldots ,m$ 
set $f_i:=c^{\mu}_w b_i$. Then multiplying the relation~(\ref{b}) 
by $c^{\mu}_w$ we get
\begin{equation}
c_w^{\lambda}c_y^{\nu}=\sum_{j=0}^m f_j,
\label{E}
\end{equation}
where the $f_j$ are weight vectors of the weight $w\lambda+y\nu$,
$f_0\in Q(y)^-$ and $ f_j\in Q(w)^+ \text { for }
j=1,\ldots,m.$ Note that $\wt_w f_i=\wt_w b_i$.

\subsubsection{}
\label{fixi}
Fix $i\,$ such that $s_i y>y .$
Then  $k[x_i] V^-_{s_i y}(\mu)=V^-_y(\mu),$ so
$f_0\in Q(y)^-\ $ implies $\ f_0.x_i^r\in  Q(s_i y)^-\ $
for any $\,r\in \Bbb {N}.$
Set $x=x_i,\ \varepsilon_i=\varepsilon .$
\subsubsection{}
\label{siwlw}
Assume that $\ s_i w<w.$ 

Since $Q(w)^+\,$ is $U_q(\frak n^+)$ invariant the relation ~(\ref{E}) implies 
that $\ (c_w^{\lambda} c_y^{\nu}).x^*\in  Q(s_i y)^- +Q(w)^+\, .$
Since $\ c_w^{\lambda} c_{s_i y}^{\nu}=(c_w^{\lambda} c_y^{\nu}).x^*\ $
up to a non-zero scalar it follows that
$\ c_w^{\lambda} c_{s_i y}^{\nu}\in  Q(s_i y)^- +Q(w)^+\, .$
\subsubsection{}
\label{siwgw}
Assume that  $\,s_i w>w.$ Then up to a non-zero scalar one has
\begin{equation}
c^{\lambda}_{s_i w} c^{\nu}_{s_i y}=(c_w^{\lambda} c_y^{\nu}).x^*=
\sum_{j=0}^m f_j.x^n, \text {\ \ for some\ \ } 
n\in\Bbb N .
\label{SW}
\end{equation}
Let us check, using~\Lem{qww}(ii), that $f_j.x^n\in Q(s_iw)^+$ 
for $j=1,\ldots,m$. Then, by~\ref{fixi}, it implies that
$$c^{\lambda}_{s_i w} c^{\nu}_{s_i y}\in Q(s_i y)^-+
Q(s_i w)^+.$$

For $a\in R^+, z=w$ or $z=s_iw$ we set $q_z(a):=\wt_z a+z^{-1}\rwt a-\lwt a$
provided the right-hand side is defined. If $q(a)$ is defined then, 
by ~\Lem{qww}(ii), $a\in Q(z)^+$ iff $q_z(a)\not=0$.
By~\ref{gneigvl} $\ \wt_{s_iw} (f_j.x^*)\,$ is defined and
\begin{equation}
\wt_w f_j-w^{-1}\rwt f_j=\wt_{s_iw} (f_j.x^*)-(s_i w)^{-1}\rwt (f_j.x^*).
\label{F}
\end{equation}
Since $\lwt (f_j.x^*)=\lwt f_j$  this implies that
$$q_{s_iw}(f_j.x^*)-q_w(f_j)=2((s_i w)^{-1}\rwt (f_j.x^*)-w^{-1}\rwt f_j).$$
Assume that $f_j.x^n=f_j.x^*$ for some $j\not=0$. Then
$$(s_i w)^{-1}\rwt (f_j.x^*)-w^{-1}\rwt f_j=(s_i w)^{-1}\rwt (c^{\lambda}_{s_i w} 
c^{\nu}_{s_i  y})-w^{-1}\rwt (c^{\lambda}_w 
c^{\nu}_ y)=0$$
so $q_{s_iw}(f_j.x^*)=q_w(f_j)$. Since $f_j\in Q(w)^+$ 
it follows that $f_j.x^*\in Q(s_iw)^+$.

Now let us show that $f_j.x^n\not=0$ iff $f_j.x^n=f_j.x^*$.
Observe that, by~\ref{aff}, the values $\wt_w f_j$ are pairwise distinct 
for $j=0,\ldots m$ so the left-hand sides of the equality ~(\ref{F}) are 
also pairwise distinct for $j=0,\ldots m$. This implies that the 
elements $\ \{f_j.x^*\}_{j=0}^m\ $
are linearly independent. 
Set $n':=\max_{0\leq j\leq m}\varepsilon (f_j)$. Then 
$$(c_w^{\lambda} c_y^{\nu}).x^{n'}=\sum_{j=0}^m f_j.x^{n'}=
\sum_{j:\varepsilon (f_j)=n'} f_j.x^*\not=0.$$
Compairing with the relation ~(\ref{SW}) we get $n'=n$ and 
$f_j.x^n\not=0$ iff $f_j.x^n=f_j.x^*$ as required.

\subsubsection{}
\label{mini6}
Set $\ s_i\star w=\max (s_i w,w).$ Since $y\leq w\,$ it follows 
(\cite{j}, A.1.7) that  $\ s_i y\leq  s_i\star w.$

Recall  our assumption that $c_w^{\lambda} c_y^{\nu}\in Q(y)^- + Q(w)^+\ $
for some pair $(y,w) : \ y\leq w.$ Suppose $y\not= w_0\,$ so there exists
$i\,$ such that $\, s_i y >y.$ Then we conclude by~\ref{siwlw},~\ref{siwgw} 
that  
$$c^{\lambda}_{s_i\star w} c^{\nu}_{s_i y}\in Q(s_i y)^- + Q(s_i\star w)^+\ ,$$
so the assumption holds for the pair $( s_i y, s_i\star w),\ $
where $y< s_i y\leq  s_i\star w.$ By induction the assumption holds for 
the pair $(w_0,w_0):\ c_{w_0}^{\lambda} c_{w_0}^{\nu}\in (Q(w_0)^-+Q(w_0)^+).$ 

However $Q(w_0)^+=(0),\ c_{w_0}\cap Q(w_0)^-=\emptyset.$ Hence
$c_{w_0}^{\lambda} c_{w_0}^{\nu}\not\in (Q(w_0)^-+Q(w_0)^+)\ $ 
which gives a contradiction.
\end{pf}

\begin{rem}{symm}
Using ~\ref{tilde} we could prove equally that the ideal 
$Q(y,w)':=\tilde Q(y,w):c_y$ is the unique minimal 
element of the component $X(y,w).$ Therefore $\,Q(y,w)=Q(y,w)'.$
\end{rem}
\subsection{}
\begin{exa}{}
The present example illustrates that in general
$$Q(s_{\alpha}, s_{\alpha}s_{\beta})\not= Q(s_{\alpha})^-+Q(s_{\alpha}s_{\beta})^+.$$

Put $\frak g=\frak {sl}_3.$ 
The diagrams below show the intersection of prime ideals 
$Q=Q(s_{\alpha})^-,\ Q(s_{\alpha}s_{\beta})^+$ of
the ring $R^+\,$ with the right modules 
$V=V^+(\omega_{\alpha}),\ V^+(\omega_{\beta}),\
V^+(\omega_{\alpha}+\omega_{\beta})=V^+(\alpha+\beta)$.

Observe that $V^+(\alpha+\beta)|_0\ $ is two dimensional.
It is spanned by a vector $c_1\,$ orthogonal to the zero weight vector
in $\ U_q(\frak b^-) u_{s_{\alpha}(\alpha+\beta)}\ $
and a vector $c_2\,$ orthogonal to the zero weight vector
in $\ U_q(\frak b^+) u_{s_{\alpha}s_{\beta}(\alpha+\beta)}\ $,
where $u_{s_{\alpha}(\alpha+\beta)}\ ,u_{s_{\alpha}s_{\beta}(\alpha+\beta)}\ $
are the extreme weight vectors of $V^+(\alpha+\beta)$ of the corresponding
weights.

In the diagram describing the pair $Q,V$ we mark with black colour
the weight vectors of $V$ belonging to $Q\cap V$. 

The ideal $\ Q(s_{\alpha})^-.$
\begin{center}
\begin{picture}(11,8)
\put(0,7){$\textstyle V^+(\omega_{\alpha})$}
\put(0,1){\circle{.2}}
\put(0,3){\circle{.2}}
\put(0,5){\circle*{.2}}

\put(.2,1.2){$\scriptstyle  c^{\omega_{\alpha}}_{w_0}$}
\put(.2,3.2){$\scriptstyle  c^{\omega_{\alpha}}_{s_{\alpha}}$}
\put(.2,5.2){$\scriptstyle c^{\omega_{\alpha}}_e$}

\put(0,1.1){\line(0,1){1.8}}
\put(0,3.1){\line(0,1){1.8}}

\put(4,7){$\textstyle V^+(\omega_{\beta})$}
\put(4,1){\circle{.2}}
\put(4,3){\circle{.2}}
\put(4,5){\circle{.2}}
\put(4.2,1.2){$\scriptstyle  c^{\omega_{\beta}}_{w_0}$}
\put(4.2,3.2){$\scriptstyle  c^{\omega_{\beta}}_{s_{\beta}}$}
\put(4.2,5.2){$\scriptstyle c^{\omega_{\beta}}_e$}

\put(4,1.1){\line(0,1){1.8}}
\put(4,3.1){\line(0,1){1.8}}

\put(9,7){$\textstyle V^+(\alpha+\beta)$}
\put(9,1){\circle{.2}}
\put(7.5,2){\circle{.2}}
\put(10.5,2.5){\circle{.2}}

\put(9.5,3.5){\circle*{.2}}
\put(8.5,3.5){\circle{.2}}
\put(7.5,5){\circle{.2}}
\put(10.5,4.5){\circle*{.2}}
\put(9,6){\circle*{.2}}

\put(9.2,1.2){$\scriptstyle  c^{\alpha+\beta}_{w_0}$}
\put(6.5,2.2){$\scriptstyle  c^{\alpha+\beta}_{s_{\beta}s_{\alpha}}$}
\put(10.7,2.7){$\scriptstyle  c^{\alpha+\beta}_{s_{\alpha}s_{\beta}}$}
\put(9.7,3.7){$\scriptstyle  c_1$}
\put(8,3.7){$\scriptstyle  c_2$}
\put(7,5.3){$\scriptstyle  c^{\alpha+\beta}_{s_{\alpha}}$}
\put(10.7,4.7){$\scriptstyle  c^{\alpha+\beta}_{s_{\beta}}$}
\put(9.2,6.2){$\scriptstyle c^{\alpha+\beta}_e$}

\put(9,1.1){\line(1,1){1.34}}
\put(9,1.1){\line(-3,2){1.34}}
\put(7.5,2.1){\line(2,3){0.86}}
\put(10.5,2.6){\line(-1,1){0.86}}
\put(9.5,3.6){\line(1,1){0.86}}
\put(8.5,3.6){\line(-2,3){0.86}}
\put(7.5,5.1){\line(3,2){1.34}}
\put(10.5,4.6){\line(-1,1){1.34}}

\end{picture}
\end{center}

The ideal $\ Q(s_{\alpha}s_{\beta})^+.$
\begin{center}
\begin{picture}(11,8)
\put(0,7){$\textstyle V^+(\omega_{\alpha})$}
\put(0,1){\circle*{.2}}
\put(0,3){\circle{.2}}
\put(0,5){\circle{.2}}
\put(.2,1.2){$\scriptstyle  c^{\omega_{\alpha}}_{w_0}$}
\put(.2,3.2){$\scriptstyle  c^{\omega_{\alpha}}_{s_{\alpha}}$}
\put(.2,5.2){$\scriptstyle c^{\omega_{\alpha}}_e$}

\put(0,1.1){\line(0,1){1.8}}
\put(0,3.1){\line(0,1){1.8}}

\put(4,7){$\textstyle V^+(\omega_{\beta})$}
\put(4,1){\circle{.2}}
\put(4,3){\circle{.2}}
\put(4,5){\circle{.2}}
\put(4.2,1.2){$\scriptstyle  c^{\omega_{\beta}}_{w_0}$}
\put(4.2,3.2){$\scriptstyle  c^{\omega_{\beta}}_{s_{\beta}}$}
\put(4.2,5.2){$\scriptstyle c^{\omega_{\beta}}_e$}

\put(4,1.1){\line(0,1){1.8}}
\put(4,3.1){\line(0,1){1.8}}

\put(9,7){$\textstyle V^+(\alpha+\beta)$}
\put(9,1){\circle*{.2}}
\put(7.5,2){\circle*{.2}}
\put(10.5,2.5){\circle{.2}}

\put(9.5,3.5){\circle{.2}}
\put(8.5,3.5){\circle*{.2}}
\put(7.5,5){\circle{.2}}
\put(10.5,4.5){\circle{.2}}
\put(9,6){\circle{.2}}

\put(9.2,1.2){$\scriptstyle  c^{\alpha+\beta}_{w_0}$}
\put(6.5,2.2){$\scriptstyle  c^{\alpha+\beta}_{s_{\beta}s_{\alpha}}$}
\put(10.7,2.7){$\scriptstyle  c^{\alpha+\beta}_{s_{\alpha}s_{\beta}}$}
\put(9.7,3.7){$\scriptstyle  c_1$}
\put(8,3.7){$\scriptstyle  c_2$}
\put(7,5.3){$\scriptstyle  c^{\alpha+\beta}_{s_{\alpha}}$}
\put(10.7,4.7){$\scriptstyle  c^{\alpha+\beta}_{s_{\beta}}$}
\put(9.2,6.2){$\scriptstyle c^{\alpha+\beta}_e$}

\put(9,1.1){\line(1,1){1.34}}
\put(9,1.1){\line(-3,2){1.34}}
\put(7.5,2.1){\line(2,3){0.86}}
\put(10.5,2.6){\line(-1,1){0.86}}
\put(9.5,3.6){\line(1,1){0.86}}
\put(8.5,3.6){\line(-2,3){0.86}}
\put(7.5,5.1){\line(3,2){1.34}}
\put(10.5,4.6){\line(-1,1){1.34}}

\end{picture}
\end{center}

Note that
$$c^{\omega_{\alpha}}_{s_{\alpha}} c^{\omega_{\beta}}_{s_{\beta}}
\in K c_1+K c_2\subset \tilde Q(s_{\alpha},s_{\alpha}s_{\beta})=
Q(s_{\alpha})^-+Q(s_{\alpha}s_{\beta})^+.$$
By~\Rem{symm} $\ Q(s_{\alpha}, s_{\alpha}s_{\beta})=
\tilde Q(s_{\alpha},s_{\alpha}s_{\beta}):c_{s_{\alpha}}\ $ so
$\,c^{\omega_{\beta}}_{s_{\beta}}\in Q(s_{\alpha}, s_{\alpha}s_{\beta}).$ 
Yet this weight vector does not belong to either $Q(s_{\alpha})^-\ $ 
nor $ Q(s_{\alpha}s_{\beta})^+\ $ and hence not to their sum. Hence
$\ Q(s_{\alpha},s_{\alpha}s_{\beta})\not=Q(s_{\alpha})^-+ Q(s_{\alpha}s_{\beta})^+$.

\end{exa}

\subsection{}
\label{minxy}
\begin{lem}{}
For all $\ (y_1,y_2)\in W\overset{w}{\diamond}W$ 
one has $\,Q(y_1,y_2)\cap c_w=\emptyset$.
\end{lem}
\begin{pf}
Suppose that $\, c_w^{\nu}\in Q(y_1,y_2)\,$ for some $\nu\in P^+(\pi)$.
This means that $c_{y_2}^{\lambda} c_w^{\nu} \in ( Q(y_1)^-+ Q(y_2)^+)\ $ 
for some $\lambda\in P^+(\pi)$.
Since $\ y_1\leq w\,$ then $\,  Q(y_1)^-\subseteq  Q(w)^-.$
Therefore $\ c_{y_2}^{\lambda} c_w^{\nu} \in ( Q(w)^-+ Q(y_2)^+)\ $ in
contradiction to ~\Prop{mini}.
\end{pf}

The ideal $Q(y_1,y_2)$ is $T$ invariant. Therefore, by~\Rem{checkq},
$$Q(y_1,y_2)_w:=\sum_{\lambda\in P^+(\pi)}c^{-\lambda}_w
( Q(y_1,y_2)\cap V^+(\lambda)) $$
is a $\check {T}$ invariant completely prime ideal of $R^w_0$ and 
$$\check {Q}(y_1,y_2)_w:=Q(y_1,y_2)_w\#\check {T}$$
is a completely prime ideal of $\check {R}^w_0$.
\subsubsection{}
\begin{cor}{minx}
(i) For each $\ (y_1,y_2)\in W\overset{w}{\diamond}W\ $ 
the component $Y_w(y_1,y_2)\,$ of $\, \Spec R_0^w\,$
has a unique minimal element $\,Q(y_1,y_2)_w$
which is a completely prime $\check {T}$ invariant ideal. 

(ii)  For each $\ (y_1,y_2)\in W\overset{w}{\diamond}W\ $ 
the component $X_w(y_1,y_2)\,$ of $\, \Spec \check{R}_0^w$
has a unique minimal element $\,\check {Q}(y_1,y_2)_w$
which is completely prime.
\end{cor}
\begin{pf}
Since $\,Q(y_1,y_2)$ is a unique minimal element of $X(y_1,y_2)$,
it follows, by~\ref{tcomponents}, that $\,Q(y_1,y_2)_w\in Y_w(y_1,y_2)$
and, moreover, it lies in all $\check T$ invariant ideals of $Y_w(y_1,y_2)$.
By~\ref{components} every $P\in Y_w(y_1,y_2)$ (resp., $P\in X_w(y_1,y_2)$)
contains some $\check T$ invariant ideal $P'\in Y_w(y_1,y_2)$.
Hence $\,Q(y_1,y_2)_w\subset P$ (resp., $\,\check {Q}(y_1,y_2)_w\subset P$)
as required. 
\end{pf}

\subsection{}
\label{pairorder}
Define an order relation on $W\diamond W$ by the formula
$$(y,z)\succeq (y',z') \text{ iff } y\leq y',\ z\geq z'.$$
The definition of $Q(y)^{\pm}\ $ implies that for $\, y\leq y'$
one has $\ Q(y)^-\subseteq Q(y')^-\,$ (resp., $\ Q(y)^+\supseteq Q(y')^+\,$).
Similarly one has

\begin{prop}{order}
(i) $\ Q(y,z)\subseteq Q(y',z')$ iff  $(y,z)\succeq (y',z')$.

(ii) $\ Q(y,z)_w\subseteq Q(y',z')_w\ $ iff  $(y,z)\succeq (y',z')$.
\end{prop}
\begin{pf}
(i) Take $\ Q(y,z)\subseteq Q(y',z')$. Then $c_{y'}^{\lambda},\
 c_{z'}^{\lambda}\not\in Q(y,z)$ for all $\lambda\in P^+(\pi)$. 
      ~\Lem{geom} implies that $\ y\leq y',\ z'\leq z$. 

Conversly, take  $\ y\leq y'$. Then 
$$ Q(y)^-+Q(z)^+\subseteq Q(y')^- +Q(z)^+\ \Rightarrow \ 
Q(y,z)=\tilde Q(y,z):c_z\subseteq \tilde Q(y',z):c_z=Q(y',z).$$
Similarly, by~\Rem{symm}, one has 
$$Q(y',z)=\tilde Q(y',z):c_{y'}\subseteq Q(y',z'):c_{y'}=Q(y',z').$$ 
Hence (i). The assertion (ii) follows from (i) and~\ref{cnn}.
\end{pf}

\subsection{}
\label{strata}
By Propositions~\ref{dcms},~\ref{mini} and Corollaries
  ~\ref{dcmrplus},~\ref{minx} we have 
the following decompositions
$$\Spec_+ R^+ =\coprod_{(y_1,y_2)\in W\diamond W} X(y_1,y_2),
\ \ \ X(y_1,y_2)^{min}=\{ Q(y_1,y_2)\},$$
$$\Spec \check {R}_0^w=\coprod_{(y_1,y_2)\in W\overset{w}{\diamond}W}
 X_w(y_1,y_2),\ \ \ 
X_w(y_1,y_2)^{min}=\{ \check {Q}(y_1,y_2)_w\},$$
$$\Spec R_0^w=\coprod_{(y_1,y_2)\in W\overset{w}{\diamond}W}
 Y_w(y_1,y_2),\ \ \ 
Y_w(y_1,y_2)^{min}=\{ Q(y_1,y_2)_w\}.$$
Let us show that the decompositions above are stratifications i.e. that
each component $X(y_1,y_2)$ (resp., $Y_w(y_1,y_2),\ X_w(y_1,y_2)$)
is locally closed and its closure $\,\overline {X}(y_1,y_2)\,$ 
(resp., $\overline{Y_w}(y_1,y_2),\ \overline{X_w}(y_1,y_2)$)
with respect to Jacobson topology is a union of components.

One has 
$$\ X(y_1,y_2)=\left\{P\in \Spec_+ R^+|\ Q(y_1,y_2)\subseteq P,
\ c_{y_1}\cap P=\emptyset,\ \ c_{y_2}\cap P=\emptyset\right\}.$$ 
Hence $\ \overline {X}(y_1,y_2)=\{P\in \Spec_+ R^+|\ Q(y_1,y_2)\subseteq P\}$
and $X(y_1,y_2)$ is locally closed.

\Prop{order}  implies that $\ X(z_1,z_2)\subseteq \overline {X}(y_1,y_2)\ $
provided  $(y_1,y_2)\succeq (z_1,z_2)$. The inverse is also true. In fact,
take $P'\in \overline {X}(y_1,y_2).$ Fix $\,(z_1,z_2)\in W\diamond W\ $ such 
that $P'\in X(z_1,z_2)$. Then  $c_{z_i}\cap Q(y_1,y_2)=\emptyset\ $
for $i=1,2$.  By ~\Lem{geom} this implies that $\ y_1\leq z_1\leq z_2\leq y_2$
that is $(y_1,y_2)\succeq (z_1,z_2)$.
The same reasoning is suitable for $\overline {X_w}(y_1,y_2),\ 
\overline {Y_w}(y_1,y_2)$.

\subsection{}
\begin{cor}{stratorder}
$$\overline {X}(y_1,y_2)=\underset{(y_1,y_2)\succeq (z_1,z_2)} 
{\coprod_{(z_1,z_2)\in W\diamond W}} X(z_1,z_2),$$
$$\overline {X_w}(y_1,y_2)=\underset{(y_1,y_2)\succeq (z_1,z_2)} 
{\coprod_{(z_1,z_2)\in W\overset{w}{\diamond}W}} X_w(z_1,z_2),$$
$$\overline {Y_w}(y_1,y_2)=\underset{(y_1,y_2)\succeq (z_1,z_2)} 
{\coprod_{(z_1,z_2)\in W\overset{w}{\diamond}W}} Y_w(z_1,z_2).$$
\end{cor}
\section{more about the strata}
All rings in this Section are noetherian. Using this and~\cite{j}, A.2.8, we will
often identify the prime spectrum of the localization $R[c^{-1}]$, $c$ being an
Ore subset of $R$, with the subset 
$$\{P\in\Spec R|\ P\cap c=\emptyset\}.$$
 
\subsection{}
In this Section we will show that the components $Y_w(y_1,y_2)$
of $\Spec R^w_0$ are isomorphic for different $w\in W$ such that
$y_1\leq w\leq y_2$.
Moreover the components $X_w(y_1,y_2)$ of $\Spec \check {R}^w_0$ 
are isomorphic to the component $X(y_1,y_2)$ of $\Spec_+ R^+$
for all $w\in W$ such that $y_1\leq w\leq y_2$. Following 
  ~\cite{j1} we identify the component $X(y_1,y_2)$ (modulo an
action of a group $\Bbb {Z}_2^l$) with the spectrum of a Laurent
polynomial ring--- see~\ref{sztwo}---~\ref{sgamma}.

All localizations considered are localizations of domains 
so the localization maps are injective. 
We will sometimes denote by the same letter an element of a ring $R$ and
its image in a localization (or in a quotient) of $R$.
\subsubsection{}
\begin{lem}{wvector}
Take $P\in X(y,w)$. Then $P\cap V^+(\nu)=Q(y,w)\cap V^+(\nu)$
for all $\nu\in P^+(\pi)$.
\end{lem}
\begin{pf}
Assume that $P\cap V^+(\nu)\not=Q(y,w)\cap V^+(\nu)$.
By~\Lem{graded} this implies that
there exists a weight vector $c^{\nu}_{\zeta}\in P\setminus Q(y,w)$.  
Choose $\lambda$ sufficiently large for $(\zeta-w\nu)$ 
(\Defn{sf}) such that $\lambda >\nu$. Then~\Lem{qww} implies that
$$c^{-\nu}_wc^{\nu}_{\zeta}=c_w^{-\lambda}c_{\xi}^{\lambda}+
c_w^{-\lambda}c_{\eta}^{\lambda},\ \text { where } \ 
c_w^{-\lambda}c_{\xi}^{\lambda}\in L(w,0),\ \ 
c_w^{-\lambda}c_{\eta}^{\lambda}\in Q(w)^+_w.$$ 
Then $c_{\eta}^{\lambda}\in Q(w)^+$ so 
$\ c_{\xi}^{\lambda}=(c^{\lambda-\nu}_wc^{\nu}_{\zeta}-c_{\eta}^{\lambda})
\in P\setminus Q(y,w)$. By~\Rem{rssaction} for all $\nu\in P(\pi)$ one has 
$\tilde{\phi}_w^{\nu}(c_w^{-\lambda}c_{\xi}^{\lambda})=
c_w^{-\lambda}c_{\xi}^{\lambda\,}$, that is
\begin{equation}
c_{\xi}^{\lambda}c_w^{\nu}=q^{(\lambda-w^{-1}\xi,\nu)}c_w^{\nu}c_{\xi}^
{\lambda}.
\label{L2}
\end{equation}
Let us show that $c_{\xi}^{\lambda}.y_{-\mu}\in P$ for all $\mu\in Q^+(\pi)$
and all elements $y_{-\mu}\in U_q(\frak b^-)$ of a weight $(-\mu)$.
We prove this by induction on $\mu\in (Q^+(\pi),\leq)$. One has
\begin{equation}
\triangle (y_{-\mu})=y_{-\mu}\otimes 1+ \tau(\mu)\otimes y_{-\mu}+
\sum_{0<\eta<\mu} k_{\eta}\tau(\eta) y_{-\mu+\eta}\otimes y_{-\eta},\ \ \ 
k_{\eta}\in K .
\label{N2}
\end{equation}
Act by $y_{-\mu}$ on the both sides of ~(\ref{L2}). Applying~(\ref{N2}) 
and induction one obtains
$$(c_{\xi}^{\lambda}.y_{-\mu})c_w^{\nu}=q^{(\lambda-w^{-1}\xi,\nu)+(\mu,w\nu)}
c_w^{\nu}(c_{\xi}^{\lambda}.y_{-\mu})\text {   mod   } P.$$
Using formula~(\ref{QW}) we get
$$(c_{\xi}^{\lambda}.y_{-\mu})c_w^{\nu}=q^{(\lambda-w^{-1}(\xi+\mu),\nu)}
c_w^{\nu}(c_{\xi}^{\lambda}.y_{-\mu})\text\  {   mod   }\  Q^+(w)\subseteq P.$$
Therefore $\ (1-q^{2(w\nu,\mu)})c_w^{\nu}(c_{\xi}^{\lambda}.y_{-\mu})\in P\ $ 
for all $\nu\in P^+(\pi)$. Hence $c_{\xi}^{\lambda}.y_{-\mu}\in P$.

Since $c_{\xi}^{\lambda}\not\in Q(y,w)$ there exists $v\in
V^-_y(\lambda)=U_q(\frak b^-)u_{y\lambda}$
such that $\xi(v)=1$. This implies that $c_y^{\lambda}=
c_{\xi}^{\lambda}.U_q(\frak b^-)$
so $c_y^{\lambda}\in P$. This contradicts $P\in X(y,w)$. 
\end{pf}
\subsubsection{}
\begin{cor}{tinv}
$$\begin{array}{l}
\text {(i)\ \ \ \   }(\Spec_+ R^+)^{T}=\left\{ Q(y,z)\right\}_{ 
(y,z)\in W\diamond W}\ .\\
\text {(ii)\ \ \  }(\Spec R^w_0)^{\check {T}}=\left\{ Q(y,z)_w\right\}_
{(y,z)\in W\overset{w}{\diamond} W}\ .\\
\text {(iii)\ \  }X_w(y,z)=\{ P\in \Spec \check {R}^w_0|\ P\cap R^w_0=Q(y,z)_w\,\}.
\end{array}$$
(iv) Take $P\in Y_w(y,z)$.
Then a weight vector $c^{-\lambda}_wc^{\lambda}_{\xi}$ belongs to $P$
iff $c^{\lambda}_{\xi}\in Q(y,z)$.
\end{cor}
\begin{pf}
The previous lemma implies (i); (ii) obtains from (i), ~\Lem{l3} 
and the diagram~(\ref{diag1}). (iii), (iv) obtain from (ii)
and~\ref{tcomponents}.
\end{pf}
\subsection{}
\label{ywz}
For any $y,w,z\in W$ let $R^{y,w,z}$ be the minimal subalgebra 
of $\Fract R^+$ containing $c_y^{-1},\ c_w^{-1},\ c_z^{-1}$. Both right and 
left action of $T$ on $R^+$ extend to $R^{y,w,z}$. Denote  the
zero component of $R^{y,w,z}$ with respect to the left $T$-action by
$R^{y,w,z}_0$.  Then the right action  of $T$ on $R^{y,w,z}_0$ extends 
to the action of $\check {T}$. Denote the corresponding
skew-product $R^{y,w,z}_0\#\check {T}$ by $\check {R}^{y,w,z}_0$.
It is clear that $\check {R}^w_0\subset \check {R}^{y,w,z}_0$.

Now take $y\leq w\leq z$. Recall that 
$$Q(y,z)_w=Q(y,z)[c_w^{-1}]\cap R^w_0,\ \ 
Q(y,z)_z=Q(y,z)[c_z^{-1}]\cap R^z_0.$$
Therefore 
$$R^{y,w,z}_0 Q(y,z)_w\supset Q(y,z)_z,\ \ 
R^{y,w,z}_0 Q(y,z)_z\supset Q(y,z)_w.$$ 
This implies that
$$R^{y,w,z}_0 Q(y,z)_w=R^{y,w,z}_0 Q(y,z)_z,\ \ 
\check {R}^{y,w,z}_0\check {Q}(y,z)_w=
\check {R}^{y,w,z}_0\check {Q}(y,z)_z.$$

For any pair $(w_1,w_2)\in W\times W$ set 
$c_{w_1,w_2}:=\{c_{w_1}^{-\lambda}c_{w_2}^{\lambda}\}_{\lambda\in P^+(\pi)}$.
\subsubsection{}
\begin{lem}{lisoxx}
Take $y\leq w\leq z$. There are canonical isomorphisms of the
Ore localizations
\begin{equation}
(R^w_0/ Q(y,z)_w)[c^{-1}_{w,z},c^{-1}_{w,y}]
\iso R^{y,w,z}_0/ (R^{y,w,z}_0 Q(y,z)_w)\iso 
(R^z_0/ Q(y,z)_z)[c^{-1}_{z,w},c^{-1}_{z,y}],
\label{iso1}
\end{equation}
\begin{equation}
(\check {R}^w_0/\check  {Q}(y,z)_w)[c^{-1}_{w,z},c^{-1}_{w,y}]
\iso \check {R}^{y,w,z}_0/ (\check {R}^{y,w,z}_0\check {Q}(y,z)_w)\iso 
(\check {R}^z_0/\check {Q}(y,z)_z)[c^{-1}_{z,w},c^{-1}_{z,y}].
\label{iso2}
\end{equation}
\end{lem}
\begin{pf}
It is sufficient to check that all the localizations are well-defined.
Observe that the image of the set $\, c_{w,z}\cup c_{w,y}\,$ 
in the quotient ring $R^w_0/ Q(y,z)_w$ consists of normal elements
so $(R^w_0/ Q(y,z)_w)[c^{-1}_{w,z},c^{-1}_{w,y}]$ is  well-defined.
  
Let us check that the image of the set $c_{z,y}\cup c_{z,w}$ 
in the quotient ring $R^z_0/ Q(y,z)_z$ is Ore. 
Since $c_w$ is Ore in $R^+$
it follows that for any $c^{\lambda}_{\xi}\in R^+, \nu\in P^+(\pi)$
there exist $c^{\mu}_{\eta}\in R^+, \nu'\in P^+(\pi)$ such that
$c^{\lambda}_{\xi}c^{\nu'}_w=c^{\nu}_wc^{\mu}_{\eta}$. By~\ref{autfi}(i)
$c^{-\lambda}_z c^{\lambda}_{\xi}$ and $c^{\lambda}_{\xi}c^{-\lambda}_z$
coincide up to a power of $q$ modulo $Q(y,z)_z$.
Therefore up to a power of $q$ one has
$$(c_z^{-\lambda}c^{\lambda}_{\xi})(c^{-\nu'}_zc^{\nu'}_w)=
c_z^{-\lambda}c^{\lambda}_{\xi}c^{\nu'}_wc^{-\nu'}_z=c_z^{-\lambda}
c^{\nu}_wc^{\mu}_{\eta}c^{-\nu'}_z=(c^{-\nu}_zc^{\nu}_w)
(c^{-\mu}_zc^{\mu}_{\eta})\ \text { mod } Q(y,z)_z.$$
Hence the image of $c_{z,w}$ is left Ore in $R^z_0/ Q(y,z)_z$. 
Similarly it is right Ore.
Since the image of the set $c_{z,y}$ in the quotient
ring $R^z_0/ Q(y,z)_z$ consists of normal elements and they commute 
up to powers of $q$ with the elements of the image of $c_{z,w}$, 
it follows that the image of $c_{z,y}\cup c_{z,w}$ 
in the quotient ring $R^z_0/ Q(y,z)_z$ is Ore. Hence
$(R^z_0/ Q(y,z)_z)[c^{-1}_{z,w},c^{-1}_{z,y}]\,$ is also well-defined.
\end{pf}
\subsubsection{}
\begin{prop}{prpisoxx}
Take $y\leq w\leq z$. 

(i) The isomorphisms~(\ref{iso1}) give rise to an order preserving bijection of 
$Y_w(y,z)$ onto $Y_z(y,z)$. 

(ii) The isomorphisms~(\ref{iso2}) give rise to an order preserving bijection of 
$X_w(y,z)$ onto $X_z(y,z)$.
\end{prop}
\begin{pf}
The definition of $Y_w(y,z)$ and  ~\Cor{minx} imply that 
$$Y_w(y,z)\cong\Spec (R^w_0/ Q(y,z)_w)[c^{-1}_{w,z},c^{-1}_{w,y}]
=\Spec R^{y,w,z}_0/ (R^{y,w,z}_0 Q(y,z)_w).$$
Taking into account ~\Lem{lisoxx} and~\Cor{tinv}(iv), we conclude that
$$\begin{array}{l}
Y_w(y,z)\cong\Spec (R^z_0/ Q(y,z)_z)[c^{-1}_{z,w},c^{-1}_{z,y}]\cong\\
\{ P\in\Spec R^z_0|\  Q(y,z)_z\subset P,\  
P\cap( c_{z,w}\cup c_{z,y})=\emptyset\}=\\
\{ P\in\Spec R^z_0|\  Q(y,z)_z\subset P,\  
P\cap c_{z,w}=\emptyset\}=Y_z(y,z)
\end{array}$$
This gives (i); the proof of (ii) is similar.
\end{pf}
\subsection{}
\label{isoxy}
\begin{prop}{prpisoxy}
For every triple $(y,w,z)$ such that $y\leq w\leq z$ 
there is an order preserving bijection of $X_w(y,z)$ onto $X(y,z)$.
\end{prop}
\begin{pf}
From the previous proposition we conclude that it is sufficient to
check the assertion for the triples $(y,z,z)$.  Fix $z\in W$.
Using notations of ~\ref{twistaut}, denote a subalgebra 
$L(z,0)\#\check {T}$ of $\check {R}^z_0$ by $\check {L}(z,0)$ 
and a subalgebra $L(z,0)\#\{c^{\nu}_w\}_{\nu\in P(\pi)}\ $
of $R^z$ by $L(z)$. Define a map $\psi:\check {L}(z,0)\to L(z)$
setting $\ \psi (a)=a,$  for $a\in L(z,0),\ \psi(\tau(\nu))=c_z^{-z^{-1}\nu}$ 
for all $\nu\in P(\pi)$.
We conclude from~\ref{twistaut} that $\psi$ is an isomorphism of algebras.
Denote by $\Psi$ the corresponding map of $\Spec\check {L}(z,0)$
onto $\Spec L(z)$.

Taking into account that $R^z Q(z)_z^+=R^z Q(z)^+$ we conclude from
Lemmas~\ref{qww},~\ref{qyw}  that
$$\check {R}^z_0=\check {Q}(z)_z^+\oplus\check {L}(z,0),
\ \ \ R^z=R^z Q(z)^+\oplus L(z).$$
Therefore there are the following bijections
$$\Psi_1:\ H_1:=\{P\in \Spec\check {R}^z_0|\ \check {Q}(z)_z^+\subset P\}
\ \to\  \Spec\check {L}(z,0),\ \ P\mapsto P\cap\check {L}(z,0),$$
with inverse $\ I\mapsto I\oplus  \check {Q}(z)_z^+ ;$ 
$$\Psi_2:\ H_2:=\{P\in \Spec R^z|\ R^zQ(z)^+\subset P\}\ \to\ 
 \Spec L(z),\ \ P\mapsto P\cap L(z),$$
with inverse $\ I\mapsto I\oplus  R^z Q(z)^+.$
Hence $(\Psi_2^{-1}\circ\Psi\circ\Psi_1)$ is a bijection of $H_1$ onto $H_2$.
Identify $X(y,z)$ and its image in $\Spec R^z$ given by the localization 
map $R^+\to R^z$. Then
$$H_1=\coprod_{y\leq z} X_z(y,z),\ \ \ H_2=\coprod_{y\leq z} X(y,z).$$
Let us show that $\,(\Psi_2^{-1}\circ\Psi\circ\Psi_1)(X_z(y,z))=X(y,z)$
for all $y\leq z$. By~\Cor{tinv}(iii) one has
$$X_z(y,z)=\{P\in \Spec\check {R}^z_0|\ P\cap R^z_0=Q(y,z)_z\}.$$
Since $Q(y,z)_z=(Q(y,z)_z\cap L(z,0))\oplus Q(z)^+$ it follows that
$$\Psi_1 (X_z(y,z))=\{P\in \Spec\check {L}(z,0)|\ 
P\cap L(z,0)=Q(y,z)_z\cap L(z,0)\}.$$
Observe that $P\cap L(z,0)=\Psi(P)\cap L(z,0)$. Therefore
$$(\Psi\circ\Psi_1) (X_z(y,z))=\{P\in \Spec L(z)|\ 
P\cap L(z,0)=Q(y,z)_z\cap L(z,0)\}.$$
Take $J\in X(y,z)$. We conclude from~\Lem{wvector}, ~\Lem{qyw} that 
$$\Psi_2(J)\cap L(z,0)=\sum_{\nu\in P^+(\pi)} c^{-\nu}_z (V^+(\nu)
\cap J)=Q(y,z)_z\cap L(z,0).$$
Hence $\Psi_2(X(y,z))\subseteq\im(\Psi\circ\Psi_1)(X_z(y,z))$. 
Since this holds for all $y\leq z$ we conclude that
$\Psi_2(X(y,z))=\im(\Psi\circ\Psi_1)(X_z(y,z))$ as required.
\end{pf}
\subsection{}
\label{calxyz}
Fix $y\leq w$. Using notations of~\ref{ywz}, denote 
$(R^w_0 /Q (y,w)_w)[c_{w,y}^{-1}]$ by $S$ and set 
$\check {S}=S\# \check {T}$. Then the canonical map $\ \check{R}_0^w\to \check {S}$
defines a bijection of $X_{w}(y,w)$ onto $\Spec \check {S}$. We calculate
$\Spec \check {S}$ in~\ref{p0}---~\ref{specd0} below.

\subsubsection{}
\label{p0}
For each $\nu\in P(\pi)$, set $z_{\nu}:=c_w^{-\nu}c_y^{\nu}\tau(y\nu+w\nu)
\in \check {S}$.
The relations ~\ref{scr2} imply that $z_{\nu}s=sz_{\nu}$ for all $s\in S$.
Since $z_{\nu}\tau(\mu)=q^{(y\nu-w\nu,\mu)}\tau(\mu)z_{\nu}$ it
follows that $z_{\nu}\in Z(\check {S})$
iff $y\nu=w\nu$. Set 
$$P_0(\pi):=\{\nu\in P(\pi)|\ y^{-1}\nu-w^{-1}\nu=0\}$$
which is a subgroup of $P(\pi)$ so that $P(\pi)/P_0(\pi)$ is 
torsion-free. Choose a subgroup $P_1(\pi)$ 
such that $P(\pi)=P_0(\pi)\oplus P_1(\pi)$. Set $T_0:=\tau (P_0(\pi)),\ 
T_1:=\tau (P_1(\pi))$. Denote the subalgebra $S\# T_0$ of
$\check {S}$ by $D$. Then $\check {S}=D\# T_1$. 

Observe that $S$ is noetherian,
so by~\cite{mr}, 2.9 $D$ is also noetherian.

\begin{lem}{restr}
The map $\psi: J\mapsto J\cap D$ is an order preserving
bijection of $\Spec \check {S}$
onto $(\Spec D)^{\check {T}}$.
\end{lem}
\begin{pf}
Since $P(\pi)=P_0(\pi)\oplus P_1(\pi)$ it follows
that $\check {T}=T_0 T_1$. Therefore $(\Spec D)^{\check {T}}=(\Spec D)^{T_1}$.
By~\Lem{l3} $\psi$ maps $\Spec \check {S}$ onto
$(\Spec D)^{T_1}$ and the map $\ I\mapsto (I\# T_1)$ is
a right inverse of $\psi$. Let us show that this is
also a left inverse of $\psi$, that is $J=(J\cap D)\# T_1$
for all $J\in\Spec \check {S}$.
Fix $J\in\Spec \check {S}$, $a\in J$. 
Write $a=\sum_{\mu}a_{\mu}\tau(\mu):\ 
\mu\in P_1(\pi),\ a_{\mu}\in D$. 
Recall that the elements $z_{\nu}$ commute with all elements
of $S$  and
$$z_{\nu}\tau(\mu)=q^{(y\nu-w\nu,\mu)}\tau(\mu)z_{\nu}=
q^{(\nu,y^{-1}\mu-w^{-1}\mu)}\tau(\mu)z_{\nu}.$$
Therefore $z_{\nu}s=sz_{\nu}$ for all $s\in D$. Since $z_{\nu}$
is invertible in $\check {S}$ one has
$$z_{\nu}az_{\nu}^{-1}=\sum_{\mu}a_{\mu}z_{\nu}\tau(\mu)z_{\nu}^{-1}=
\sum_{\mu}q^{(\nu,y^{-1}\mu-w^{-1}\mu)}a_{\mu}\tau(\mu)\in J.$$
The values $(y^{-1}\mu-w^{-1}\mu)$ are pairwise distinct
for different $\mu\in P_1(\pi)$, so $a_{\mu}\tau(\mu)\in J$
for all $\mu\in P_1(\pi)$. Then $a_{\mu}\in J\cap D$ and
$J=(J\cap D)\# T_1$ as required.
\end{pf}
\subsubsection{}
\label{sztwo}
Let $r$ be the rank of $P_0(\pi)$. Identify $P_0(\pi)/2P_0(\pi)$
with $\Bbb {Z}_2^r$. For each $\tau(\nu)\in T_0$ let
$d(\tau(\nu))$ denote the image of $\nu$ in $\Bbb {Z}_2^r$.
For $s\in S$ set $d(s):=0$. This defines $\Bbb {Z}_2^r$
grading on $D$. For $g\in \Bbb {Z}_2^r$ denote the subspace
$\{a\in D|\ d(a)=g\}$ by $D_g$. Denote by $\Gamma$ the character
group of $\Bbb {Z}_2^r$. For each $\gamma\in\Gamma$ define
$\theta_{\gamma}\in\Aut D$ setting 
$\theta_{\gamma}|_{D_g}:=\gamma(g)\cdot\id$. 
View $\Gamma$ as acting on ideals of $D$ via the $\theta_{\gamma}:\ 
\gamma\in\Gamma$ and hence on $\Spec D$. Since the $\theta_{\gamma}$
commute with the action of $\check {T}$ it follows that $\Gamma$
acts also on $(\Spec D)^{\check {T}}$.

\begin{lem}{ztwo}
The map taking $I\in (\Spec D_0)^{\check {T}}$ to the minimal
primes over $DI$ (with inverse $P\mapsto P\cap D_0$)
is a bijection of $(\Spec D_0)^{\check {T}}$
onto the $\Gamma$ orbits of $(\Spec D)^{\check {T}}$. 
\end{lem}
\begin{pf}
Since $D=S\# T_0$ it follows that $D=D_0 T_0=T_0 D_0$.
This implies that $DI$ is a two-sided graded ideal of $D$ 
for any $T$ invariant ideal $I$ of $D_0$. The reasoning
of~\cite{j}, 1.3.9 implies that for any 
$I\in (\Spec D_0)^{\check {T}}$ the minimal
primes $Q_i$ over $DI$ form a single $\Gamma$ orbit and satisfy
$I=Q_i\cap D_0$ for all $i$. 

Let us show that the inverse map is well-defined.
Fix $P\in (\Spec D)^{\check {T}}$ and set
$I:=P\cap D_0$. Assume that $I$ is not prime. Then, by~\Lem{l2},
there exist homogeneous $a,b\in D_0\setminus I$ such that
$aD_0b\subseteq I$. Then $\ aDb=aD_0T_0b=aD_0bT_0\subseteq IT_0\subseteq P\,$
that contradicts $P$ being prime and completes the proof.
\end{pf}

\begin{rem}{rankz2}
For $i=1,\ldots,l$ define the element $\sigma_i\in\Aut\check {R}_0^w$
by the formulas
$$ \sigma_i|_{R_0^w}=\id;\ \sigma_i(\tau(\omega_i))=-\tau(\omega_i);\
\sigma_i(\tau(\omega_j))=\tau(\omega_j)\text{ for } j\not=i.$$

Consider the group ${\Bbb Z}_2^l\subseteq\Aut\check {R}_0^w$ generated
by the automorphisms $\sigma_i$. This group acts naturally on $D$ and
the image of ${\Bbb Z}_2^l$ in $\Aut D$ identifies with $\Gamma$.
\end{rem}

\subsubsection{}
\label{specd0}
Denote the subalgebra of $\check {S}$ generated
by the central elements $z_{\nu}:y\nu=w\nu$ by $Z$.
Take $\mu\in P(\pi)$; then $\mu=y\nu+w\nu$ for some $\nu$ such that
$y\nu=w\nu$ iff $\mu\in 2P_0(\pi)$.
It follows that $Z\subset D_0$ and $Z$ is a Laurent polynomial
ring of the rank $r$.
Since $D_0=S\# \tau(2P_0(\pi))$ it follows that $D_0\cong S\otimes Z$
as $\check {T}$ algebras (the action of $\check {T}$ on $Z$ is trivial).
Since $S$ is noetherian, $D_0$ is also noetherian.

\begin{lem}{d0}
(i) The map $P\mapsto P\cap Z$ is an isomorphism of $(\Spec D_0)^{\check {T}}$
onto $\Spec Z$.

(ii) For each $P\in (\Spec D_0)^{\check {T}}$, the quotient
$D_0/P$ is a domain.
\end{lem}
\begin{pf}
Take $P\in (\Spec D_0)^{\check {T}}$.
Since $P$ is prime and $Z$ is contained in the centre of $D_0$
one has $(P\cap Z)\in\Spec Z$.

Take any $I\in\Spec Z$. Since $Z$ is $\check {T}$ invariant
then $Q:=SI$ a two-sided  $\check {T}$ invariant ideal
of $D_0$ contained in $P$. Identify $D_0$ with  $S\otimes Z$.
Then $Q=S\otimes I$ and $D_0/Q\cong S\otimes (Z/I)$ as $\check {T}$ 
algebras, where the action of $\check {T}$ on $Z/I$ is trivial. 
Since $Z/I$ is a domain, $G:=(Z/I)\setminus\{ 0\}$ is an Ore subset of 
$S\otimes (Z/I)$. Set $F:=\Fract (Z/I)$ and identify $S\otimes (Z/I)[G^{-1}]$
with $S\otimes F$. The action of $\check {T}$ on $S\otimes (Z/I)$ extends
to $S\otimes F$. 
By definition $S=(R^w_0 /Q (y,w)_w)[c_{w,y}^{-1}]$. This is
a domain for any choice of the base field
$K\supseteq k(q)$. Set (for a moment) $K:=F$. Then we get that
$S\otimes F$ is a domain so $S\otimes (Z/I)$ is also a domain.
Hence $Q$ is a completely prime ideal of $D_0$. Since $Q\cap Z=I$
this establishes the surjectivity in (i).

Take $P\in (\Spec D_0)^{\check {T}}$ and set $I:=(P\cap Z)$.
Again set $Q=SI$ and define $G,F$ as above.
Denote by $\overline P$ the image of $P/Q$ in $S\otimes (Z/I)$ which
is a prime $\check {T}$ invariant ideal. Recall that $P\cap Z=I$ so 
$\overline P\cap G=\emptyset$. Hence $\overline P[G^{-1}]$ is a prime 
$\check {T}$ invariant ideal of $S\otimes (Z/I)[G^{-1}]=S\otimes F$.
\Cor{tinv}(ii) implies that the zero ideal is the only $\check {T}$ invariant 
prime ideal of the ring $S\otimes K'$ for any field $K'$
containing $k(q)$. Hence $\overline P=(0)$ that is $P=Q$. 
This establishes (ii) and injectivity in (i).
\end{pf}

\subsubsection{}
\label{sgamma}
Recall that $Z$ is a subalgebra of $(\check {R}^w_0 /Q (y,w)_w)[c_{w,y}^{-1}]$
generated by the central elements
$z_{\nu}:=c_w^{-\nu}c_y^{\nu}\tau(2w\nu)$ where $\nu\in P(\pi)$
such that $y\nu=w\nu$. For each $z\in W$ denote by $r(z)$
the rank of the free group $P_z(\pi):=\{\mu\in P(\pi)|\ z\mu=\mu\}$
(one has $r(z)=l-s(z)$, where $s(z)$ denotes the minimal length of 
an expression for $z$ as a product of reflections).
Then $\rk Z=r(w^{-1}y)$. 
Combining~\ref{calxyz}---~\ref{specd0} one obtains the 

\begin{prop}{prpgamma}
The map $P\mapsto (P/Q (y,w)_w)[c_{w,y}^{-1}]\cap Z$ is an isomorphism
of the space of ${\Bbb Z}_2^l$ orbits in $X_w(y,w)$ onto $\Spec Z$.
\end{prop}

Now Propositions~\ref{prpisoxx},~\ref{prpisoxy},~\ref{prpgamma} give the

\begin{thm}{spectra}
$$\text {(i) }\Spec_+ R^+ =\coprod_{(y,z)\in W\diamond W} X(y,z), $$
where each $X(y,z)$ is isomorphic up to an action of $\Bbb {Z}_2^l$
to the spectrum of the Laurent polynomail ring of rank $r(y^{-1}z)$.
$$\text {(ii) }\Spec \check{R}_0^w=\coprod_{(y,z)\in W\overset{w}{\diamond}W}
X_w(y,z),$$
where each $X_w(y,z)$ is isomorphic to the component $X(y,z)$ of 
$\Spec_+ R^+ $.
\end{thm}

\section{The Centre of $R_0^w$}
\label{ccc}
Denote the element $(c^{\lambda}_{\xi})^{-1}$ of  $\Fract R^+$ by 
$c^{-\lambda}_{\xi}$. Set 
$$A:=\{a\in \Fract R^+|\ c^{\lambda}_{\xi}a\in R^+ 
\text { for some } \lambda\in P^+(\pi),\ \xi\in \Omega(V(\lambda)^*)\}.$$
The right action of $U_q$ on $R^+$ extends to $A$ and $a=c^{-\lambda}_{\xi}b$
is a weight vector iff $b\in R^+$ is a weight vector.
\subsection{}

\begin{lem}{weightcentre}
Let $a$ be a weight vector of $A$.
Then $a\in Z(\Fract R_0^e)$ iff $a\in K c^{-\nu}_{e}c^{\nu}_{w_0}$
for some $\nu\in P(\pi)$ satisfying $w_0\nu=-\nu$.
\end{lem}
\begin{pf}
By~\cite{j}, 9.1.4(i), 10.1.11(ii) for any $\nu, \lambda,\in P^+(\pi),
 \mu\in\Omega(V^+(\lambda))$ one has 
\begin{equation}
\ c^{\lambda}_{\mu} c^{\nu}_e=
q^{(\nu ,\mu- \lambda )} c^{\nu}_e c^{\lambda}_{\mu} ,\ \ \ \  
\ c^{\lambda}_{\mu}c^{\nu}_{w_0}=q^{-(w_0\nu ,\mu-w_0\lambda )} 
c^{\nu}_{w_0}c^{\lambda}_{\mu}.
\label{N}
\end{equation}
This implies that $\ \ 
c^{w_0\nu}_{w_0}c^{\nu}_eb=bc^{w_0\nu}_{w_0}c^{\nu}_e\ \ $ 
for any $\nu\in P(\pi),\ b\in R_0^e$. 
Hence  $\ c^{-\nu}_{e}c^{\nu}_{w_0}\in Z(\Fract R_0^e)\ $ if 
$w_0\nu=-\nu$.

Let us prove the converse. 
For each $b\in A$ consider the set of pairs $\ \{ (\lambda,\xi)\in P^+(\pi)\times 
\Omega(V^+(\lambda))|\ c^{\lambda}_{\xi}b\in R^+\}$. This set admits a lexicographic
preorder $(\lambda,\xi)\leq (\lambda',\xi')$ iff $\lambda\leq\lambda'$ or
$\lambda=\lambda'\ \text{and } \xi\leq\xi'$. The expression
 $b=c^{-\lambda}_{\xi} d$ ($d\in R^+$) will be called {\em a reduced decomposition} 
if the pair $(\lambda,\xi)$ is a minimal with respect to the preorder above.

Set 
$$B:=\{b\in R^+|\ b\not\in c_{w_0}^{\omega_i} R^+,\ 
b\not\in c_e^{\omega_i} R^+\ \text{ for all }i=1,\ldots ,l\}.$$
Given $b\in R^+$ write $b=c^{\nu_1}_{w_0}c^{\nu_2}_eb':
\ \nu_1,\nu_2\in P^+(\pi),\ b'\in B$. Theorem 3 of~\cite{j3} implies that
$Q(w_0s_i)^+=c_{w_0}^{\omega_i} R^+\,$ 
(similarly $Q(s_i)^-=c_e^{\omega_i} R^+\,$). Since $Q(w_0s_i)^+,Q(s_i)^-$
are completely prime ideals of $R^+$ it follows that $\nu_1,\nu_2, b'$ are
uniquely determined. The element $b'$ will be called {\em the abnormal part} of $b$.

Let $a$ be a non-zero weight vector of $A$ and let $a\in Z(\Fract R_0^e)$.
Fix a reduced decomposition $a=c^{-\lambda}_{\xi}d$. Let $c^{\lambda_1}_{\mu_1},\ 
c^{\lambda_2}_{\mu_2}\,$ be the abnormal parts of $c^{\lambda}_{\xi},d$ 
respectively. One has 
$$a=c^{-\lambda}_{\xi}d=q^rc^{\nu_1}_{w_0}c^{\nu_2}_ec^{-\lambda_1}_{\mu_1}
c^{\lambda_2}_{\mu_2}\ \ \ \text { for some  } \nu_1,\nu_2 \in P(\pi), r\in \Bbb Z.$$
Set $b:=c^{-\lambda_1}_{\mu_1}c^{\lambda_2}_{\mu_2}$. Observe that 
$b=c^{-\lambda_1}_{\mu_1}c^{\lambda_2}_{\mu_2}$ is a reduced
decomposition.

Let $c^{\nu}_{\eta}$ be a weight vector of $R^+$. One has
$\ c_e^{-\nu}c^{\nu}_{\eta}a=ac_e^{-\nu}c^{\nu}_{\eta}.$

The relations ~(\ref{N}) imply that 
\begin{equation}
c^{\nu}_{\eta}b=q^r bc^{\nu}_{\eta} \text {\ \ \ \    for some   } r\in\Bbb Z. 
\label{K}
\end{equation}
Moreover one has
$$b c^{\omega_i}_e=q^{(\wt_e b,\omega_i)}c^{\omega_i}_e b, \text { where } \wt_e
b=\mu_2-\mu_1-\lambda_2+\lambda_1.$$
Act by $x_i$ on the both sides of the relation above. Taking into account that 
$\wt_e (b.x_i)=\wt_e b-\alpha_i$ we obtain
$$q^{(\wt_e b,\omega_i)}(1-q^{-2})c^{\omega_i}_e (b.x_i)=b c^{\omega_i}_{s_i}-q^{(\wt_e
b,\omega_i)-(\alpha_i,\rwt b)}c^{\omega_i}_{s_i}b.$$
Using~(\ref{K}) we conclude that $c^{\omega_i}_e (b.x_i)\in K b c^{\omega_i}_{s_i}$.
One has 
$$b.x_i=(c^{-\lambda_1}_{\mu_1}c^{\lambda_2}_{\mu_2}).x_i=
c^{-\lambda_1}_{\mu_1}(c^{\lambda_2}_{\mu_2}.x_i)
-q^{(\alpha_i,\mu_1-\mu_2)}
c^{-\lambda_1}_{\mu_1}(c^{\lambda_1}_{\mu_1}.x_i)(c^{-\lambda_1}_{\mu_1}
c^{\lambda_2}_{\mu_2})=
(c^{\lambda_1}_{\mu_1})^{-2} d \text { for some } d\in R^+.$$
Therefore 
$$c^{\omega_i}_e (c^{\lambda_1}_{\mu_1})^{-2} d\in K b c^{\omega_i}_{s_i} \ \ 
\Rightarrow \ \ 
c^{\omega_i}_e d\in K  (c^{\lambda_1}_{\mu_1})^2 b c^{\omega_i}_{s_i}=
K c^{\lambda_1}_{\mu_1}c^{\lambda_2}_{\mu_2} c^{\omega_i}_{s_i}.$$

Recall that  $c^{\lambda_1}_{\mu_1}, c^{\lambda_2}_{\mu_2}\in B$ so
$\ c^{\lambda_1}_{\mu_1}c^{\lambda_2}_{\mu_2} c^{\omega_i}_{s_i}\not\in Q(s_i)^-$. 
Since $c^{\omega_i}_e\in Q(s_i)^-$ it follows that $d=0$ so $b.x_i=0$. 
Replacing $ c^{\omega_i}_e$ by  $c^{-w_0\omega_i}_{w_0}$
and $Q(s_i)^-$ by $Q(s_i w_0)^+$ we get $b.y_i=0$. Since $b.x_i=b.y_i=0$
it follows that $b.t_i=b$. 

Let us check that $\lambda_1=0$.Assume the converse.
Then $c^{\lambda_1}_{\mu_1}\not=c^{\lambda_1}_{w_0}$ 
since $c^{\lambda_1}_{\mu_1}\in B$. Therefore there exists $i$ such that 
$c^{\lambda_1}_{\mu_1}.x_i\not =0$. 

Since $b.x_i=0,\ b.t_i=b$ one has
$$c^{\lambda_2}_{\mu_2}.x_i=(c^{\lambda_1}_{\mu_1}b).x_i=
(c^{\lambda_1}_{\mu_1}.x_i)b\ \ \Rightarrow
\ \ b=(c^{\lambda_1}_{\mu_1}.x_i)^{-1}(c^{\lambda_2}_{\mu_2}.x_i).$$

Yet $\rwt (c^{\lambda_1}_{\mu_1}.x_i)<\rwt c^{\lambda_1}_{\mu_1}\ $ this contradicts
$b=c^{-\lambda_1}_{\mu_1}c^{\lambda_2}_{\mu_2}$ being a reduced decomposition.

Now, $\lambda_1=0$ and therefore $c^{\lambda_2}_{\mu_2}.x_i=
c^{\lambda_2}_{\mu_2}.y_i=0$ 
for all $i=1,\ldots , l$. Then $\lambda_2=0$ so $b\in K^*$. Hence 
$a\in K^*c^{\nu_1}_{w_0}c^{\nu_2}_e$. Since $a\in Z(\Fract R_0^e)\,$
it follows that $\nu_1+\nu_2=0$. Moreover the relations ~(\ref{N}) 
imply that $\ (\nu_2,\mu)-
(w_0\nu_1,\mu)=0\,$ for any $\mu\in Q^-(\pi)$. Hence 
$a\in K^*c^{\nu_1}_{w_0}c^{-\nu_1}_e$ and $\nu_1+w_0\nu_1=0$ as required.
\end{pf}

\subsection{}
\label{theta}
Let $\theta$ be the automorphism of the Dynkin diagram 
defined by the property
$w_0 \omega_i=-\omega_{\theta(i)}$. One has $\theta^2=1$.
Set 
$$\frak I:=\left\{ i\in \left\{1,\ldots ,l\right\}|\theta(i)=i\right\},\ 
\overline {\frak I}:=\left\{ i\in \left\{1,\ldots ,l\right\}|\theta(i)>i\right\}.$$
Set $\ z_i:=c_e^{-\omega_i}c_{w_0}^{\omega_i}$; for $i\in \overline {\frak I}$ set
$\tilde z_i:=z_i z_{\theta(i)}.$ 

One has $z_i\in R^e_0,\ z_i^{-1}\in R_0^{w_0}$.
For $w=e$ the centre $Z(R_0^e)$ is the polynomial algebra generated by the set
$\, M:=\left\{z_i:\ i\in\frak I,\ \tilde z_i :\ i\in \overline {\frak
I}\right\}$ --- see~\cite{j}, 7.1.20. Similarly $Z(R_0^{w_0})$ is the 
polynomial algebra
generated by the set $M^{-1}=\left\{ m^{-1}:\ m\in M\right\}$.
We will show that $Z(R_0^w)$ is the polynomial algebra generated by the set
$(M\cup M^{-1})\cap R_0^w$.

For a more precise description of the set of generators of $Z(R_0^w)$ set 
$${\frak I}_w^-:=\left\{i\in \frak I|\ w\omega_i=\omega_i\right\},\ 
\overline {\frak I}_w^-:=\left\{i\in \overline{\frak I}|\ w\omega_i=\omega_i,\ 
 w\omega_{\theta (i)}=\omega_{\theta (i)}\right\},\ $$

$${\frak I}_w^+:=\left\{i\in \frak I|\ w\omega_i=w_0\omega_i\right\},\ 
\overline {\frak I}_w^+:=\left\{i\in \overline{\frak I}|\ w\omega_i=w_0\omega_i,\ 
 w\omega_{\theta (i)}=w_0\omega_{\theta (i)}\right\}.$$

Then $M\cap R_0^w=\left\{z_i:\ i\in\frak I_w^-,\ \tilde z_i :\ i\in \overline {\frak
I}_w^-\right\}$ and $M^{-1}\cap R_0^w=
\left\{z_i^{-1}:\ i\in\frak I^+_w,\ \tilde z_i^{-1}:\ i\in \overline {\frak
I}^+_w\right\}$. 
\subsubsection{}
\begin{prop}{centreis}
The centre $Z(R_0^w)$ is the polynomial algebra generated by the set
$C:=(M\cup M^{-1})\cap R_0^w$.
\end{prop}
\begin{pf}
Set $z_{\nu}:=c^{-\nu}_ec^{\nu}_{w_0}$ for all $\nu\in P(\pi)$ satisfying $w_0\nu=-\nu$.
Observe that $R^w_0\subset A$. Then, 
in view of ~\Lem{weightcentre}, it suffices to show that any element 
$z_{\nu}\in R^w_0$ can be expressed as a product of elements of $C$.

Write $\nu=\sum k_i\omega_i$ and set $\frak A^-:=\{i:k_i<0\},\ \frak
A^+:=\{i:k_i>0\}$. 
Set 
$$\nu_1:=-\sum_{i\in\frak A^-} k_i\omega_i,\ \ \ \ \ \ 
\nu_2:=\sum_{i\in\frak A^+} k_i\omega_i.$$
Then $\nu=\nu_2-\nu_1,\ \nu_1,\nu_2\in P^+(\pi)$. Since $w_0\nu=-\nu$ it 
follows that $k_{\theta(i)}=k_i$ so $\theta (\frak A^{\pm})=\frak A^{\pm}$
and $w_0\nu_1=-\nu_1,\ w_0\nu_2=-\nu_2$. Hence
$$z_{\nu}=z_{\nu_1}^{-1}z_{\nu_2};\ \ \ \ \ \ \ \ 
z_{\nu_1}=\prod_{i\in\frak I\cap\frak A^-}z_i^{k_i}
\prod_{i\in\overline{\frak I}\cap\frak A^-}\tilde {z}_i^{k_i},\ \ \ \ \ 
z_{\nu_2}=\prod_{i\in\frak I\cap\frak A^+}z_i^{k_i}
\prod_{i\in\overline{\frak I}\cap\frak A^+}\tilde {z}_i^{k_i}.$$
Let us show that $z_i\in R^w_0$ for all $i\in \frak A^+$(then also
$\tilde {z}_i\in R^w_0$ for $i\in\frak I\cap\frak A^+$) and
$z_i^{-1}\in R^w_0$ for all $i\in \frak A^-$.

Observe that $z_i\in R^w_0$ if $w\omega_i=\omega_i$ and 
$z_i^{-1}\in R^w_0$ if $w\omega_i=w_0\omega_i$. Hence it suffices
to check that $w\omega_i=\omega_i$ (resp. $w\omega_i=w_0\omega_i$) for all 
$i\in \frak A^+$ (resp. $i\in \frak A^-$).

Since $z_{\nu}\in R^w_0$ there exists $\lambda\in P^+(\pi)$ such that
$$c^{-\lambda}_wc^{\lambda}_{\xi}=z_{\nu}=c^{-\nu_1}_{w_0} c^{\nu_2}_{w_0} 
c^{\nu_1}_e c^{-\nu_2}_e\ \  \Rightarrow
\ \ c^{\lambda}_{\xi}c^{\nu_1}_{w_0}c^{\nu_2}_e\in K^* 
c^{\lambda}_w c^{\nu_2}_{w_0}c^{\nu_1}_e.$$
Take $i\in \frak A^+$. Then $c^{\nu_2}_e\in Q(s_i)^-$, whereas 
$c^{\nu_2}_{w_0}c^{\nu_1}_e\not\in Q(s_i)^-$. From the formula above
we conclude that $c^{\lambda}_w\in Q(s_i)^-$ so $w\omega_i=\omega_i$. 
Similarly $i\in \frak A^-$ implies that $w\omega_i=w_0\omega_i$.
\end{pf}
\subsubsection{}
\begin{rem}{nonisomorphic}
\Prop{centreis} implies that the rings $R^w_0$ are in general non-isomorphic:
they have centres of different Gelfand-Kirillov dimension. Observe that 
this dimension is maximal if $w=e,w_0$. If $\frak g$ is simple then
for all $w\not=e,w_0$ one has  $\dim Z(R^w_0)<\dim Z(R^e_0)$.

In fact, fix $w$ is such that 
$\dim Z(R^w_0)=\dim Z(R^e_0)$. This implies that
$\overline {\frak I}=\overline {\frak I}_w^-\cup\overline {\frak I}_w^+$
and $\frak I=\frak {I}_w^-\cup\frak {I}_w^+$. Set $\frak {J}_1:=
\{i|\ w\omega_i=\omega_i\},\ \frak {J}_2:=
\{i|\ w\omega_i=w_0\omega_i\}$. Then $\frak {J}_1\cup\frak {J}_2=\{1,\ldots,l\},\ 
\frak {J}_1\cap\frak {J}_2=\emptyset$.
Observe that
$w\in W_2$ where $W_2$ is a subgroup of $W$ which is generated  by
$\{s_i:\ i|\ w\omega_i\not=\omega_i \}=\{s_i:\ i\in \frak {J}_2\}$.
Similarly, $w_0w\in W_1$ where $W_1$ is a subgroup of $W$ which is generated  by
$\{s_i:\ i|\ (w_0w)\omega_i\not=\omega_i \}=\{s_i:\ i\in \frak {J}_1\}$.
Since $w_0=(w_0w)w^{-1}$ it follows that $w_0\in W_1W_2$ so $W=W_1W_2$.
Since $\frak g$ is simple, one has either $W=W_1$ or $W=W_2$. This means that
$w=e$ or $w=w_0$.
\end{rem}

\section{appendix: index of notations}
Symbols used frequently are given below under the section number
where they are first defined.

\ref{ch}\ \ \ \ \ \ $k,K,U_q(\frak g),\check {T},\check {U}_q(\frak g),U_q(\frak n^-),
x_i,y_i,t_i^{\pm 1},l,W,S^w$

\ref{dcpm}\ \ \ \ \ \ $w_0$

\ref{second}\ \ \ \ \ \ $\pi,Q(\pi),Q^{\pm}(\pi),\omega_i,P(\pi),\geq,P^+(\pi),\tau,
V(\lambda), c^{\lambda}_{\xi,v},$
$$\ \ R_q[G], V^+(\lambda),R^+,
\Omega ( V^+(\lambda)),c^{\lambda}_w,c^{\lambda}_{\xi},c_w,R^w,R^w_0,\check {R}^w_0$$

\ref{normal}\ \ \ \ \ \ $A[c^{-1}]$

\ref{defskewpr}\ \ \ \ \ \ $\# $

\ref{scr}\ \ \ \ \ \ \ \ $\lwt ,\rwt ,\cdot |_{\lambda},\ \cdot |^{\lambda},
\ c^{\lambda}_{\mu}$

\ref{scr1}\ \ \ \ \ \ $J^{\pm}_{\lambda}(\eta)$

\ref{scr1w}\ \ \ \ \ \ $J^{\pm}_{\lambda}(\eta)_w$

\ref{autfi}\ \ \ \ \ \ $\phi_w^{\nu},\Phi_w$

\ref{fiinv}\ \ \ \ \ \ $D^{\pm}_P(\nu)$

\ref{specplus}\ \ \ \ \ \ $P^{++},\Spec_+ R^+$
 
\ref{dcm11}\ \ \ \ $X(y_-,y_+)$

\ref{tspect}\ \ \ \ $(\Spec_+ R^+)^T,\ (\Spec R^w)^T$

\ref{demaz}\ \ \ \ $V_y^{\pm}(\lambda),V_y^{\pm}(\lambda)^{\perp},Q(y)^{\pm}$

\ref{dcm1}\ \ \ \ $W\diamond W$

\ref{cnn}\ \ \ \ $\Spec_w R^+$

\ref{components}\ \ \ \ $W\overset{w}{\diamond}W,X_w(y_1,y_2),Y_w(y_1,y_2)$

\ref{epsilon}\ \ \ \ $U,\varphi _i,\varepsilon _i,y_i^*,x_i^*,y_w^*,x_w^*$

\ref{ntt}\ \ \ \ $Q(y)^{\pm}_w$

\ref{defwtw}\ \ \ \ \ \ $\wt_w$

\ref{qywy}\ \ \ \ \ \ $Q(y,w)_w$

\ref{qinr}\ \ \ \ \ \ $Q(y,w)$

\ref{minxy}\ \ \ \ \ $Q(y_1,y_2)_w,\check {Q}(y_1,y_2)_w$

\ref{pairorder}\ \ \ \ \ $\succeq $ 

\ref{ywz}\ \ \ \ \ \ $c_{w_1,w_2}$

\ref{calxyz}\ \ \ \ \ \ $S,\check {S}$

\ref{p0}\ \ \ \ $z_{\nu},P_0(\pi),P_1(\pi),T_0,T_1,D$

\ref{sztwo}\ \ \ \ $D_0,\Gamma$

\ref{specd0}\ \ \ \ $Z$



\end{document}